\documentclass[10pt, reqno]{amsart}
\usepackage{amssymb,latexsym,amsfonts}
\usepackage{amsthm}
\usepackage{amsmath,calligra,mathrsfs}
\usepackage{color}
\usepackage[OT2, T1]{fontenc}
\usepackage{url}
\usepackage{array}
\usepackage{hyperref}
\usepackage{paralist}
\usepackage{mathabx}			% for \widec
\usepackage{amscd}
\usepackage[all,cmtip]{xy}			% for commutative diagrams
\usepackage{sseq}		
\usepackage{verbatim}
\usepackage{parskip}
\usepackage{microtype}
\usepackage[in]{fullpage}
\usepackage{mathrsfs} 			% for \mathscr (script letters)
\usepackage{cleveref}
\usepackage{enumitem}
\usepackage{moreenum}			% for enumeration via Greek letters
\usepackage{colonequals}	
\usepackage{fontenc}
\usepackage{amsmath}
\usepackage{geometry}
\usepackage{graphicx}
\usepackage{stmaryrd}
\usepackage[all]{xy}
\usepackage[french]{babel}
\usepackage[latin1]{inputenc}
\usepackage[T1]{fontenc}
\usepackage[foot]{amsaddr}
\usepackage{subfiles}
\usepackage{ragged2e}
\usepackage[bbgreekl]{mathbbol}     % This package is in order to have the prism symbol
\DeclareSymbolFontAlphabet{\mathbb}{AMSb}   % This is in order to common blackboard letters in the usual style (else the prism package messes them up)

\DeclareMathOperator{\Mod}{Mod}

\DeclareMathOperator{\bD}{\mathbb{D}}

\DeclareMathOperator{\cL}{\mathcal{L}}

\DeclareMathOperator{\Spec}{Spec}

\DeclareMathOperator{\inv}{inv}
\DeclareMathOperator{\nil}{nil}

\DeclareMathOperator{\Hom}{Hom}

\DeclareMathOperator{\End}{End}

\DeclareMathOperator{\Par}{Par}

\DeclareMathOperator{\codim}{codim}

\DeclareMathOperator{\Gr}{Gr}
\DeclareMathOperator{\Spc}{Spc}
\DeclareMathOperator{\Aff}{Aff}
\DeclareMathOperator{\Sh}{Sh}

\DeclareMathOperator{\obl}{obl}
\DeclareMathOperator{\ind}{ind}

\DeclareMathOperator{\Vect}{Vect}
\DeclareMathOperator{\Rhom}{R\mathscr{H}\text{\kern -3pt {\calligra\large om}}\,}
\DeclareMathOperator{\EExt}{\mathscr{E}\text{\kern -3pt {\calligra\large xt}}\,}
\DeclareMathOperator{\Ind}{Ind}
\DeclareMathOperator{\Aut}{Aut}

\DeclareMathOperator{\ad}{ad}
\DeclareMathOperator{\Ad}{Ad}
\DeclareMathOperator{\Ens}{Ens}

\DeclareMathOperator{\FLoc}{FLoc}

\DeclareMathOperator{\sinf}{<}

\DeclareMathOperator{\Perv}{Perv}
\DeclareMathOperator{\Loc}{Loc}

\DeclareMathOperator{\Ker}{Ker}

\DeclareMathOperator{\pt}{pt}

\DeclareMathOperator{\colim}{colim}
\DeclareMathOperator{\Fl}{Fl}

\DeclareMathOperator{\Res}{Res}
\DeclareMathOperator{\res}{res}

\DeclareMathOperator{\bqlk}{\overline{\mathbb{Q}}_{\ell,\kappa}}

\DeclareMathOperator{\coinv}{coinv}
\DeclareMathOperator{\St}{St}
\DeclareMathOperator{\Rep}{Rep}
\DeclareMathOperator{\PrSh}{PrSh}
\DeclareMathOperator{\PreSh}{PrSh}
\DeclareMathOperator{\AlgSp}{AlgSp}

\DeclareMathOperator{\Kos}{Kos}

\DeclareMathOperator{\PrStk}{PrStk}

\DeclareMathOperator{\IndLoc}{IndLoc}

\DeclareMathOperator{\triv}{triv}

\DeclareMathOperator{\AlgAss}{AlgAss}

\DeclareMathOperator{\gr}{gr}

\DeclareMathOperator{\cT}{\mathcal{T}}

\DeclareMathOperator{\act}{act}

\DeclareMathOperator{\Lie}{Lie}

\DeclareMathOperator{\Fonct}{Fonct}
\DeclareMathOperator{\Frac}{Frac}

\DeclareMathOperator{\PrCat}{PrCat_{\ell}}
\DeclareMathOperator{\Catst}{Cat_{st,\ell}}

\DeclareMathOperator{\clp}{\mathcal{L}^{+} }

\usepackage{aliascnt}
\newaliascnt{numberingbase}{subsubsection}
\numberwithin{equation}{numberingbase}

\newtheoremstyle{thms}{0.7em}{0pt}{\itshape}{}{\bfseries}{.}{ }{}
\theoremstyle{thms}
\newtheorem{conj}[numberingbase]{Conjecture}
\newtheorem{cor}[numberingbase]{Corollaire}
\newtheorem{lemma}[numberingbase]{Lemme}

\newtheorem{prop}[numberingbase]{Proposition}
\newtheorem{Q}[numberingbase]{Question}
\newtheorem{thm}[numberingbase]{Théorème}

\newtheoremstyle{claims}{0.7em}{0pt}{}{}{\itshape}{.}{ }{}
\theoremstyle{claims}
\newtheorem{claim}[equation]{Claim}

\newtheoremstyle{defs}{0.7em}{0pt}{}{}{\bfseries}{.}{ }{}
\theoremstyle{defs}
\newtheorem{defi}[numberingbase]{Définition}

\newtheorem{exa}[numberingbase]{Exemple}
\newtheorem*{exas}{Exemples}
\newtheorem{rmq}[numberingbase]{Remarque}
\newtheorem*{rmqs}{Remarques}

\Crefname{claim}{Claim}{Claims}
\Crefname{conj}{Conjecture}{Conjectures}
\Crefname{cor}{Corollaire}{Corollaires}
\Crefname{defn}{Définition}{Définitions}
\Crefname{eg}{Exemple}{Exemples}
\Crefname{prop}{Proposition}{Propositions} 
\Crefname{Q}{Question}{Questions}
\Crefname{rem}{Remarques}{Remarques}
\Crefname{theorem}{Théorème}{Théorèmes}
\Crefname{thm}{Théorème}{Théorèmes}
\Crefname{variant}{Variant}{Variants}

\theoremstyle{thms}
\newtheorem{thm-tweak}[subsection]{Théorème}
\Crefname{thm-tweak}{théorème}{théorèmes}
\newtheorem{lemma-tweak}[subsection]{Lemme}
\Crefname{lemma-tweak}{lemme}{Lemmes}
\newtheorem{cor-tweak}[subsection]{Corollaire}
\Crefname{cor-tweak}{corollaire}{corollaires}
\newtheorem{prop-tweak}[subsection]{Proposition}
\Crefname{prop-tweak}{proposition}{propositions} 
\newtheorem{conj-tweak}[subsection]{Conjecture}
\Crefname{conj-tweak}{Conjecture}{Conjectures} 

\theoremstyle{defs}
\newtheorem{defn-tweak}[subsection]{Définition}
\Crefname{defn-tweak}{définition}{définitions}
\newtheorem{eg-tweak}[subsection]{Exemple}
\Crefname{eg-tweak}{exemple}{exemples}
\newtheorem*{rmqs-tweak}{Remarques}
\newtheorem{rmq-tweak}[subsection]{Remarque}
\Crefname{rmq-tweak}{remarque}{remarques}

\newtheoremstyle{subsection-tweak}
   {11pt}
   {3pt}%
   {}
   {}%
   {\bfseries}
   {}%
   {.5em}
   {\thmnumber{\@{#1}{}\@{#2}.}%
    \thmnote{~{\bfseries#3.}}}

\theoremstyle{subsection-tweak}
\newtheorem{pp}[numberingbase]{}
\newcommand{\bpp}{\begin{pp}}
\newcommand{\epp}{\end{pp}}

\theoremstyle{subsection-tweak}
\newtheorem{pp-tweak}[subsection]{}

\renewcommand{\b}{\textbf}

\newcommand{\brems}{\begin{rmqs} \hfill \begin{enumerate}[label=\b{\thenumberingbase.},ref=\thenumberingbase]}
\newcommand{\remi}{\addtocounter{numberingbase}{1} \item}
\newcommand{\erems}{\end{enumerate} \end{rmqs}}
\newcommand{\bexas}{\begin{exas} \hfill \begin{enumerate}[label=\b{\thenumberingbase.},ref=\thenumberingbase]}

\newcommand{\eexas}{\end{enumerate} \end{exas}}

\newcommand{\bsm}{\begin{smallmatrix}}
\newcommand{\esm}{\end{smallmatrix}}
\newcommand{\blem}{\begin{lemma}}
\newcommand{\elem}{\end{lemma}}

\newcommand{\bconj}{\begin{conj}}
\newcommand{\econj}{\end{conj}}
\newcommand{\bprob}{\begin{Problem}}
\newcommand{\eprob}{\end{Problem}}

\newcommand{\bq}{\begin{Q}}
\newcommand{\eq}{\end{Q}}
\newcommand{\benum}{\begin{enumerate}[label={{\upshape(\alph*)}}]}
\newcommand{\benuma}{\begin{enumerate}[label={{\upshape(\arabic*)}}]}
\newcommand{\benumr}{\begin{enumerate}[label={{\upshape(\roman*)}}]}
\newcommand{\eenum}{\end{enumerate}}
\newcommand{\bitem}{\begin{itemize}}
\newcommand{\eitem}{\end{itemize}}
\newcommand{\bc}{\begin{comment}}
\newcommand{\ec}{\end{comment}}
\newcommand{\bdefi}{\begin{defi}}
\newcommand{\edefi}{\end{defi}}
\newcommand{\ov}{\overline}

\newcommand{\bexa}{\begin{exa}}
\newcommand{\eexa}{\end{exa}}

\newcommand{\bcl}{\begin{claim}}
\newcommand{\ecl}{\end{claim}}

\newcommand{\ba}{\begin{aligned}}
\newcommand{\ea}{\end{aligned}}
\newcommand{\be}{\begin{equation}}
\newcommand{\ee}{\end{equation}}
\newcommand{\ti}{\widetilde}
\newcommand{\bpf}{\begin{proof}}
\newcommand{\epf}{\end{proof}}
\newcommand{\bthm}{\begin{thm}}
\newcommand{\ethm}{\end{thm}}
\newcommand{\bprop}{\begin{prop}}
\newcommand{\eprop}{\end{prop}}
\newcommand{\bcor}{\begin{cor}}
\newcommand{\ecor}{\end{cor}}
\newcommand{\brem}{\begin{rmq}}
\newcommand{\erem}{\end{rmq}}
\newcommand{\Gm}{\Gamma}

\def\dar[#1]{\ar@<1pt>[#1]\ar@<-2pt>[#1]}
  \entrymodifiers={!!<0pt,0.7ex>+}
\makeatletter
\def\@tocline#1#2#3#4#5#6#7{
    \begingroup 
    \@ifempty{#4}{}{}

    \parindent\z@ \leftskip#3\relax \advance\leftskip\@tempdima\relax
    #5\hskip-\@tempdima
      \ifcase #1
       \or\or \hskip 2em \or \hskip 1em \else \hskip 3em \fi%
      #6\nobreak\relax
    \dotfill\hbox to\@pnumwidth{\@tocpagenum{#7}}\par
    \nobreak
    \endgroup
 }
 \def\l@section{\@tocline{1}{0pt}{1pc}{}{}}

\renewcommand{\tocsection}[3]{%
  \indentlabel{\@ifnotempty{#2}{\makebox[1.3em][l]{%
    \ignorespaces#1 \bfseries{#2}.\hfill}}}\bfseries{#3}
    \vspace{1.5pt}}

\renewcommand{\tocsubsection}[3]{%
  \indentlabel{\@ifnotempty{#2}{\hspace*{-0.5em}\makebox[2.1em][l]{%
    \ignorespaces#1#2.\hfill}}}#3
    \vspace{1.5pt}}
	
\makeatother 

\makeatletter % This is for proper formatting of appendices in the table of contents
\newcommand\appendix@section[1]{%
  \refstepcounter{section}%
  \orig@section*{Appendix \@Alph\c@section. #1}%
}
\let\orig@section\section
\g@addto@macro\appendix{\let\section\appendix@section}
\makeatother

\newcommand{\ra}{\rightarrow}
\newcommand{\hra}{\hookrightarrow}
\newcommand{\cH}{\mathcal{H}}

\newcommand{\cU}{\mathcal{U}}
\newcommand{\cC}{\mathcal{C}}

\newcommand{\cM}{\mathcal{M}}
\newcommand{\cY}{\mathcal{Y}}

\newcommand{\co}{\mathcal{O}}

\newcommand{\cO}{\mathcal{O}}
\newcommand{\cS}{\mathcal{S}}
\newcommand{\cA}{\mathcal{A}}

\newcommand{\g}{\gamma}
\newcommand{\la}{\lambda}

\newcommand{\kt}{\mathfrak{t}}
\newcommand{\kb}{\mathfrak{b}}

\newcommand{\eps}{\epsilon}

\newcommand{\cI}{\mathcal{I}}

\newcommand{\kg}{\mathfrak{g}}

\newcommand{\kD}{\mathfrak{D}}

\newcommand{\bZ}{\mathbb{Z}}

\newcommand{\bql}{\overline{\mathbb{Q}}_{\ell}}

\newcommand{\kc}{\mathfrak{c}}

\newcommand{\kC}{\mathfrak{C}}
\newcommand{\wkC}{\widetilde{\mathfrak{C}}}
\newcommand{\ab}{\mathbb{A}}
\newcommand{\ev}{ev}

\newcommand{\NN}{\mathbb{N}}
\newcommand{\al}{\alpha}
\newcommand{\bG}{\mathbb{G}}
\newcommand{\J}{\mathrm{J}}

\newcommand{\br}{\textbf{r}}

\newcommand{\cF}{\mathcal{F}}
\newcommand{\cD}{\mathcal{D}}
\newcommand{\cP}{\mathcal{P}}
\newcommand{\cX}{\mathcal{X}}
\newcommand{\La}{\Lambda}
\title{Faisceaux caractères sur les espaces de lacets d'algèbres de Lie}
\author{Alexis Bouthier}
\begin{document}
\maketitle
\begin{center}
\textbf{Abstract:}
\end{center}
 We establish several foundational results regarding the Grothendieck-Springer affine fibration. More precisely, we prove some constructibility results on the affine Grothendieck-Springer sheaf and its coinvariants, enrich it with a group of symmetries, analog to the situation of Hitchin fibration, prove some perversity statements once we take some derived coinvariants and construct some specialization morphisms for the homology of affine Springer fibers.
\section{Introduction}
\subsection{Préambule}
\subsubsection{Rappels sur les résultats de \cite{BKV}}
Soit un corps $k$ algébriquement clos, $G$ un groupe connexe réductif sur $k$, $(B,T)$ une paire de Borel et $W$ le groupe de Weyl associé.
La théorie de Springer usuelle s'attache à construire des faisceaux pervers $\Ad(G)$-équivariants sur $\Lie(G)$ ainsi que des représentations de $W$. Pour ce faire, on s'intéresse à la fibration de Grothendieck-Springer:
\[\pi:\tilde{\kg}=\{(g,\g)\in G/B\times\kg, \ad(g)^{-1}\g\in\Lie(B)\}\ra\kg,\]
 qui est petite, projective et  finie étale de groupe $W$ au-dessus du lieu régulier semisimple $j:\kg^{rs}\hra\kg$. En considérant le faisceau de Grothendieck-Springer $\cS_{fin}=\pi_{*}\bql[\dim(\tilde{\kg})]$, on obtient donc un faisceau pervers qui est l'extension intermédiaire de sa restriction à $\kg^{rs}$ et tel que $\End(\cS_{fin})=\bql[W]$. En particulier, $\cS_{fin}$ est muni d'une action de $W$.
Pour chaque représentation irréductible $V$ de $W$, on obtient alors une composante $V$-isotypique $\cS_{fin,V}$ qui est également perverse  et $\Ad(G)$-équivariante. Ces faisceaux sont les analogues pour l'algèbre de Lie des faisceaux caractères de Lusztig. Plus généralement, pour obtenir suffisamment de faisceaux caractères, il est important de considérer des faisceaux $\cS_{fin,\cL}=\pi_{*}\cL[\dim(\tilde{\kg})]$ où $\cL$ est un système local.

Dans \cite{BKV}, Bouthier, Kazhdan et Varshasky entreprennent l'extension de cette théorie dans le cas affine, i.e. $G$ est remplacé par le groupe de lacets $G((t))$, $\kg$ par $\kg((t))$ et $W$ par le groupe de Weyl affine étendu $\widetilde{W}=X_*(T)\rtimes W$, où $X_{*}(T)$ désigne les cocaractères du tore maximal $T$.
La fibration de Grothendieck-Springer est remplacée par son analogue affine:
\[f:\wkC=\{(g,\g)\in G((t))/I\times\kg((t)),\ad(g)^{-1}\g\in\Lie(I)\}\ra\kg((t))\] 
où $I$ est l'Iwahori associé à $B$. En fait, la fibration se factorise par le sous-ind-schéma $\kC$ des éléments compacts.

 Dans ce contexte, tous les objets considérés sont des ind-schémas, il est donc plus délicat de donner un sens à la notion de \og petitesse\fg, de construire une théorie des faisceaux ou même de parler de perversité. 
A ceci s'ajoute que les  fibres géométriques réduites sont des $k$-schémas localement de type fini dont la cohomologie, bien que d'amplitude bornée, n'est pas de dimension finie.
Dans \cite[Thm. 7.1.4]{BKV}, on montre qu'en un certain sens, au-dessus du lieu génériquement régulier semisimple $\kC_{\bullet}$, $f$ est un morphisme petit en un sens raisonnable et que si l'on travaille avec des faisceaux $G((t))$-équivariants, on dispose de foncteurs cohomologiques avec une notion de perversité. On considère alors  la flèche induite:
\[[f]:[\wkC_{\bullet}/G((t))]\ra[\kC_{\bullet}/G((t))],\]
et le complexe $\cS=[f]_!\omega_{\wkC_{\bullet}}$ où $\omega_{\wkC_{\bullet}}$ est le dualisant. On montre que ce complexe est pervers, s'obtient comme extension intermédiaire de sa restriction à un lieu régulier semisimple $\kC_{rs}$ et  vérifie que $\End(\cS)=\bql[\widetilde{W}]$.
\subsubsection{Les prolongements de \cite{BKV}}
Plusieurs questions et applications sont laissées en jachère dans \cite{BKV}, la première de celles-ci, ainsi que mentionnée dans \cite[0.5.1]{BKV} est le passage aux coinvariants dérivés. Cette deuxième étape est fondamentale puisqu'elle réalise le projet de construire de faisceaux caractères qui géométrisent certains caractères de représentations de groupes $p$-adiques.
Pour chaque représentation de dimension finie $V$ de $\widetilde{W}$, on considère la composante $V$-isotypique dérivée et on s'attend à obtenir un faisceau pervers constructible. 
Etant donné que $\widetilde{W}$ est maintenant infini, le foncteur des coinvariants n'est plus exact, et il n'est pas clair que l'on obtienne un faisceau pervers.
De plus, pour espérer obtenir un complexe constructible, des énoncés de finitude du faisceau $\cS$ vu comme complexe de $\bql[\widetilde{W}]$-modules sont nécessaires.
Enfin, de même que dans le cas classique, il est important de considérer des faisceaux plus généraux que $\omega_{\wkC_{\bullet}}$ pour lesquels on a des énoncés de perversité.
Le présent travail s'attache donc à établir ces résultats et améliore les attentes de \cite{BKV} où l'on pensait qu'il était nécessaire d'avoir la pureté des fibres de Springer affines pour obtenir la perversité des coinvariants \cite[0.5.1]{BKV}.
\subsection{Symétries locales et perversité}
Pour être en mesure de pouvoir calculer les coinvariants dérivés sous $\widetilde{W}$, on a besoin de connaître la façon dont $\widetilde{W}$ agit. Par analogie avec la théorie de Springer globale, initiée par Yun \cite{Yun} qui concerne la fibration de Hitchin parabolique $f^{par}:\cM^{par}\ra\cA^{par}$, la partie sphérique, i.e. l'action de la sous-algèbre $\bql[X_{*}(T)]^{W}$ de $\bql[\widetilde{W}]$ se factorise via l'action du faisceau des composantes connexes $\pi_{0}(\cP^{glob}/\cA^{par})$ du champ de Picard $\cP^{glob}$ qui agit sur $\cM^{par}$ au-dessus de $\cA^{par}$. De plus, d'après Lusztig \cite{Lu}, on a a priori une action différente de $\widetilde{W}$ sur l'homologie des fibres de Springer affines que celle construite dans \cite{BKV} et il est important de pouvoir comparer les deux.
On construit donc un Picard local $\cP$ qui agit sur $\wkC$ au-dessus de $\kC$ (\ref{P-act},\ref{P-act2}) et on a le théorème suivant \ref{lusz},\ref{P-act3}:
\bthm
Le faisceau de Grothendieck-Springer affine $\cS$ est naturellement $\cP$-équivariant et cette action commute à celle de $\widetilde{W}$. De plus, l'action de $\widetilde{W}$ construite par extension intermédiaire coïncide avec celle de Lusztig.
\ethm
On en profite également pour généraliser l'énoncé de perversité de \cite{BKV} pour des systèmes locaux.
Plus spécifiquement, on a une équivalence de faisceaux étales $[\wkC/G((t))]\cong[\Lie(I)/I]$. Le champ $\cX=[\Lie(I)/I]$ s'écrit alors comme limite projective de champs d'Artin lisses de type fini $\cX \simeq\varprojlim\cX_i$ avec des flèches de transition lisses et on considère une catégorie des systèmes locaux renormalisés :
\[\Loc^{ren}(\cX)\simeq\colim_{f^{!}}\Loc^{ren}(\cX_i),\]
avec $\Loc^{ren}(\cX_i)=\{\cL\otimes_{\bql}\omega_{\cX_i},\cL\in\Loc(Y)\}$ (on montre que cette définition ne dépend pas des choix).
On a alors le théorème suivant \ref{fond-spr2}:
\bthm
Pour tout $\cL\in\Loc^{ren}([\wkC_{\bullet}/G((t))])$, $[f]_{!}\cL$ est pervers et s'obtient comme extension intermédiaire de sa restriction à $\kC_{rs}$.
\ethm
\subsection{Constructibilité}
Comme on l'a déjà vu, l'homologie des fibres de Springer affines bien que d'amplitude finie, est de dimension infinie en chaque degré en général. Toutefois, une fois que l'on prend les coinvariants dérivés sous $\widetilde{W}$, on s'attend à obtenir des faisceaux pervers constructibles, dont la trace de Frobenius doit nous donner des caractères de représentations. On montre en fait un énoncé plus fort de finitude avant de prendre les coinvariants \ref{W-constr}:
\bthm
Le faisceau de Grothendieck-Springer affine $\cS$ est constructible comme faisceau de $\bql[\widetilde{W}]$-modules.
En particulier, pour toute représentation de dimension finie $\tau$ de $\widetilde{W}$, le faisceau $\cS_{\tau}$ des $\tau$-coinvariants est un faisceau constructible sur $[\kC_{\bullet}/G((t))]$ (cf. définition \eqref{tfdef1}.)
\ethm
Indépendamment des applications de cet énoncé au présent article, il s'avère également utile pour des applications fines sur le calcul du support singulier du faisceau de Grothendieck-Springer affine et son lien avec la conjecture de locale constance de Goresky-Kottwitz-McPherson \cite[sect. 4.4]{Bt}.

Maintenant que l'on a un énoncé de finitude pour les $\tau$-coinvariants, on s'intéresse à la perversité de ceux-ci.
\subsection{Lemme d'homotopie et perversité des coinvariants}
\subsubsection{La situation fibre à fibre de Yun}
Dans le cas de la fibration de Hitchin globale, on a par le lemme d'homotopie \cite[Lem. 3.2.3]{LN} que l'action de $\cP^{glob}$ sur les faisceaux de cohomologie $R^{i}f_{*}^{par}\bql$ se factorise par $\pi_{0}(\cP^{glob})$. De plus, Yun définit un morphisme de faisceaux :
\[\sigma:\bql[X_*(T)]^{W}\ra\pi_{0}(\cP^{glob})\]
et il montre en particulier (\cite[Thms. 1, 2 ,3]{Yun}) que la partie sphérique de l'action de $\widetilde{W}$ sur la cohomologie à support propre $R^{i}f_{*}^{par}\bql$ des fibres de Hitchin se factorise par $\sigma$.
Il montre également un énoncé local, mais qui vaut fibre à fibre. A savoir pout tout $\g\in\kC_{\bullet}(k)$, on a un morphisme local:
\[\sigma_{\g}:\bql[X_*(T)]^{W}\ra\pi_{0}(\cP_{\g})\]
et à nouveau la partie sphérique agit sur $H_c^{i}(f^{-1}(\g),\bql)$ via $\sigma_{\g}$. En ce qui concerne l'homologie, l'énoncé est  plus faible.
\subsubsection{L'énoncé $\infty$-catégorique au niveau des complexes}
Dans le cas qui nous concerne, on a besoin d'une généralisation de ce résultat de Yun dans plusieurs directions. Tout d'abord, pour calculer les coinvariants dérivés, il faut travailler au niveau des $\infty$-catégories, d'une part, et d'autre part, on a besoin d'un lemme d'homotopie qui porte sur le complexe plutôt que sur ses faisceaux de cohomologie.
En effet, pour le moment, il est délicat de définir, s'il existe, un faisceau local $\pi_{0}(\cP)$ des composantes connexes, ce qui nous force à travailler fibre à fibre et donc on ne peut utiliser la perversité de $\cS$.
On montre donc  un énoncé d'homotopie général  \ref{cor-homo} pour l'action d'un schéma en groupe lisse commutatif sur un complexe.
\bthm
Soit un $k$-schéma de type fini $S$, $P$ un $S$-schéma en groupes lisse. On considère $K\in\cD(S)$ un faisceau $P$-équivariant, alors l'action se factorise par le faisceau des composantes connexes relatives $\pi_{0}(P/S)$, i.e. pour tout $\phi:U\ra S$ l'action de $P(U)$ sur $\phi^!K$ se factorise par $\pi_0(P/S)(U)$.
\ethm
Bien que la preuve de cet énoncé soit relativement élémentaire, il est à noter qu'à notre connaissance, aucun énoncé analogue n'était connu pour l'action d'un schéma en groupes sur des complexes. Enfin, contrairement à ce qui était attendu (\cite[0.5.1]{BKV}), cet énoncé permet de contourner la pureté des fibres de Springer affines, pour en déduire la perversité des coinvariants \ref{perv-co}.

On tire ensuite de cet énoncé général le corollaire suivant le corollaire qui s'applique aussi bien au cas global que local  \ref{Pic-pi}, \ref{fact}:
\bcor
\benum
\item
L'action de $\cP^{glob}$ sur le complexe $f^{par}_*\bql$ se factorise en une action de $\pi_{0}(\cP^{glob}/\cA)$.
\item
Pour tout point géométrique $\g\in\kC_{\bullet}(k)$, l'action de $\cP_{\g}$ sur  $R\Gm(f^{-1}(\g),\omega_{X_{f^{-1}(\g)}})$ se factorise par $\pi_0(\cP_{\g})$.
\eenum
\ecor
On remarquera que contrairement à ce qui était attendu, il n'est pas nécessaire de connaître la pureté des fibres de Springer affines pour obtenir cet énoncé d'homotopie. 

Cela nous permet également de formuler, sous forme conjecturale, un énoncé de Yun généralisé:
\bconj\label{conjI}
Pour tout corps algébriquement clos $K$ et $\g\in\kC_{\bullet}(K)$, l'action de la partie sphérique $\bql[X_*(T)]^W$ de $\widetilde{W}$ sur $R\Gm(f^{-1}(\g),\omega_{f^{-1}(\g)})$ se factorise par \eqref{sigma}.
\econj
Notons qu'un certain nombre d'énoncés globaux établis par Yun devraient s'étendre tels quels, au niveau des complexes et $\infty$-catégorique, le point délicat de cette conjecture est le passage du local au global et la preuve de cette conjecture bien qu'abordable nécessite un travail spécifique.

En supposant cette conjecture, on va montrer un énoncé de perversité  pour les coinvariants. Si $\tau$ est une représentation de dimension finie de $\widetilde{W}$, on dit qu'elle est de torsion si $\tau_{\vert X_*(T)}$ est une somme directe de caractères de torsion.
On montre alors \ref{coinv}:
\bthm\label{perv-co}
Soit $G$ semisimple simplement connexe, soit $\tau$ une représentation de torsion de dimension finie de $\widetilde{W}$, on suppose \ref{conjI} vérifiée, alors le faisceau des $\tau$-coinvariants $\cS_{\tau}$ est pervers.
\ethm
\subsubsection{Morphismes de spécialisation}
Enfin dans la dernière section, on construit des morphismes de spécialisation pour l'homologie des fibres de Springer affines compatibles à l'action de $\widetilde{W}$ (\ref{sp-sprin}). Ce type de morphisme permet en particulier de relier l'homologie des fibres de Springer affines aux éléments déployés à celle pour des éléments elliptiques et devrait s'avérer utile pour étudier des familles de fibres de Springer affines, de la même manière que Ngô l'a fait avec les fibres de Hitchin.
\subsubsection{Remerciements}
Ce travail s'inscrit dans la continuation de ma collaboration avec David Kazhdan et Yakov Varshavsky. J'ai en particulier beaucoup bénéficié de nos discussions ainsi que de leurs suggestions et idées. Je remercie également pour nos multiples entretiens Gérard Laumon et Eric Vasserot.
Enfin, je remercie Sam Raskin, Marco Robalo et  Christophe Cornut d'avoir aimablement répondu à mes questions.
\section{Fibration de Grothendieck-Springer affine}

\subsection{Quotient adjoint et centralisateur régulier}
On commence par faire quelques rappels sur le centralisateur régulier. Les hypothèses que l'on fait dans la suite ne sont pas optimales, aussi bien sur la caractéristique que sur la base, mais suffisantes pour ce travail.
Les énoncés sont tirés de Ngô \cite[1.1-1.5]{N} et Bouthier-\v{C}esnavi\v{c}ius \cite[sect.4]{BC}, pour les énoncés les plus généraux.

\subsubsection{Théorème de Chevalley}
Soit un corps $k$ algébriquement clos, un groupe $G$ connexe réductif sur $k$, $\kg$ son algèbre de Lie, $(B,T)$ une paire de Borel, $W$ le groupe de Weyl et $\kt=\Lie(T)$. On suppose que la caractéristique est première à l'ordre de $W$.
La restriction des fonctions induit un isomorphisme de Chevalley:
\[k[\kg]^{G}\stackrel{\sim}{\rightarrow}k[\kt]^{W}.\]
On en déduit un morphisme polynôme caractéristique \[
\chi:\kg\ra\kc=\kt/W\]
 ainsi qu'un morphisme $\pi:\kt\ra\kt/W$. L'espace $\kc$ est lisse affine. Soit $\kD=\prod_{\al\in R} d\al\in k[\kt]^{W}$ le diviseur discriminant où $R$ désigne l'ensemble des racines. On note $\kc^{rs}$ le complémentaire du discriminant ainsi que $\kt^{rs}$ et $\kg^{rs}$ leurs images respectives dans $\kt$ et $\kg$ via les morphismes $\chi$ et $\pi$.
Au-dessus de $\kc^{rs}$, $\pi$ est fini étale de groupe $W$.
Soit le schéma des centralisateurs $C=\{(g,\g)\in G\times\kg,\ad(g)\g=\g\}$, on définit l'ouvert régulier $\kg^{reg}=\{(g,\g)\in G\times\kg, \dim C_{\g}=r\}$.
Au-dessus de cet ouvert régulier la flèche:
\[\chi_{reg}:\kg^{reg}\ra\kc\]
est lisse et les fibres géométriques sont des espaces homogènes sous $G$.
On dispose d'une section à $\chi_{reg}$ dite de Kostant:
\[\eps:\kc\ra\kg^{reg},\]
qui généralise la matrice compagnon pour les groupes classiques.
\subsubsection{Centralisateur régulier}
Au-dessus de l'ouvert régulier, le schéma des centralisateurs $C_{\vert\kg^{reg}}$ est lisse et se descend en un schéma en groupes lisse commutatif $J\ra\kc$ de telle sorte que l'on a un isomorphisme canonique:
\[\chi_{reg}^{*}J\stackrel{\sim}{\rightarrow} C_{\vert\kg^{reg}}.\]
qui s'étend en un morphisme:
\begin{equation}
j_{\kg}:\chi^{*}J\rightarrow C.
\label{carJ}
\end{equation}
De manière explicite, $J=\eps^{*}C$  pour $\eps$ une section de Kostant.
En particulier, la flèche:
\begin{equation}
G\times\kg^{reg}\ra\kg^{reg}\times_{\kc}\kg^{reg}
\label{Jtors}
\end{equation}
donnée par  $(g,\g)\mapsto (\g,\ad(g).\g)$ est un $J$-torseur.
On a besoin d'une version qui prend en compte le Borel. La proposition suivante est due à Yun \cite[Prop. 2.3.1]{Yun}:
\bprop\label{B-var}
Soit $J_{\kb}$ le tiré-en-arrière à $\kb=\Lie(B)$, soit $C_{\kb}$ le schéma des centralisateurs au-dessus de $\kb$ pour l'action adjointe de $B$.
Alors on a un morphisme canonique:
\[j_{\kb}:J_{\kb}\ra C_{\kb}\]
de telle sorte que $j_{\kg,\vert\kb}$ se factorise en :
\[J_{\kb}\stackrel{j_{\kb}}{\rightarrow}C_{\kb}\ra C_{\vert \kb}.\]
\eprop
\subsection{Objets placides}
\subsubsection{Espaces algébriques placides}
\bdefi\label{fpl}
\benumr
\remi
 Soit un morphisme d'espaces algébriques $f:X\ra S$, il admet une présentation placide ou est placidement présenté (resp. fortement pro-lisse) s'il existe une présentation filtrante $X\simeq\varprojlim X_{\al}$, où les morphismes de transition sont affines lisses et les $X_{\al}$ sont des $S$-espaces algébriques de présentation finie (resp. et de plus les $X_{\al}$ sont lisses sur $S$).
\remi
Un morphisme  d'espaces algébriques $X\ra Y$ est dit placide s'il existe $Y'\ra Y$ étale surjectif, tel que $X'=X\times_{Y}Y'\simeq\coprod_{\al}X_{\al}$ où $X_{\al}\ra Y$ est placidement présenté.
\remi
Si $S$ est le spectre d'un corps, on dit que $X$ est placide (resp. admet une présentation placide, resp. fortement pro-lisse) si la flèche $X\ra\Spec(k)$ l'est.
\remi
Un morphisme de schémas $X\ra Y$ est dit lisse, si localement pour la topologie étale, on a $X=\coprod X_{\al}$ avec $X_{\al}\ra Y$ fortement pro-lisse.
\eenum
\edefi
 \brems\label{fpl2}
\remi
D'après \cite[1.1.3]{BKV}, la classe des morphismes fortement pro-lisses est stable par changement de base, composition, morphismes fortement pro-lisses.
\remi
D'après \cite[1.1.6]{BKV}, si $f:X\ra S$ est placidement présenté, alors on a une présentation canonique $X\cong\varprojlim_{X\ra Y} Y$, où la limite parcourt les morphismes fortement pro-lisses $X\ra Y$ pour $Y$ de présentation finie sur $S$ et les flèches de transition sont lisses affines.
\erems

\bexa\label{arc}
Soit un anneau commutatif $A$ et $X$ un $A[[t]]$-schéma affine de type fini.
Pour tout entier $n\in\NN$, on considère le foncteur des arcs tronqués
\[\cL_n X:B\mapsto X(B[t]/(t^{n+1})).\]
Il est représentable par un $A$-schéma affine de type fini et on forme l'espace d'arcs $\clp X=\varprojlim\cL_n X$.
Il représente le foncteur $B\mapsto X(B[[t]])$ et est un $A$-schéma affine non-nécessairement noethérien.
On dispose alors d'une flèche de projection $\clp X\ra X$.
\eexa
\subsubsection{Ind-schémas ind-placides}\label{ind-sc-ind-pl}
Soit  un anneau commutatif $A$, on appelle ind-schéma sur $A$, tout foncteur sur les $A$-schémas affines $T:(\Aff_A)^{op}\ra\Ens$ qui est représentable par une colimite filtrante de $A$-schémas qcqs $T_{\al}$, où les morphismes de transition sont des immersions fermées.
D'après \cite[Lem. 1.1.4]{Gai}, si $\cC_{/T}$ désigne la catégorie des schémas affines $S$ munis d'une immersion fermée $S\hra T$, la catégorie est filtrante et l'on a une écriture canonique pour $T$:
\[T=\colim_{S\in\cC_{/T}}S.\]

Un ind-schéma $T$ est  \textsl{ind-placide} s'il s'écrit  $T\simeq\colim T_{\al}$ avec des flèches de transition qui sont des immersions fermées de présentation finie et les $T_{\al}$ sont des $k$-schémas placides qcqs. 
L'exemple de base d'un tel ind-schéma est le suivant:
\bexa\label{loop}
Si $X$ est un $A((t))$-schéma affine de type fini, on considère le foncteur de lacets:
\[\cL X:B\mapsto X(B((t))).\]
Ce foncteur commute aux produits et pour $X=\ab^{1}$, $\cL\ab^{1}\cong\colim (\clp\ab^{1}\stackrel{t}{\rightarrow}\clp\ab^{1}\dots)$.
Enfin, si $X\hra Y$ est une immersion fermée, alors $\cL X\hra \cL Y$ est aussi une immersion fermée de telle sorte qu'en plongeant $X$ dans un $\ab^N$, $\cL X$ est représentable par un ind-schéma ind-affine.
De plus, si $A=k$ est un corps et $X$ lisse sur $k((t))$ alors $\cL X$ est un ind-schéma ind-placide $\cL X$ (\cite[Thm. 6.3]{Dr}).
\eexa
De même, on a une notion d'ind-espace algébrique et de ind-placidité pour ceux-ci.

Enfin, soit $\cX$ un champ, i.e. un faisceau étale en groupoïdes, on dit qu'il est \textsl{placide}, s'il admet un atlas $f:X\ra\cX$ où $X$ est un espace algébrique placide et $f$ une flèche représentable, surjective et lisse au sens de \ref{fpl}.
\subsection{$\infty$-champs}
\subsubsection{Généralités}
Soit une $\infty$-catégorie $C$, $\Spc$ l'$\infty$-catégorie  des $\infty$-groupoïdes. On  considère la catégorie des préfaisceaux  $\PreSh(C)=\Fonct(C^{op},\Spc)$. Si $C$ est équipée d'une topologie de Grothendieck $\cT$, on considère $\Sh(C)$ la sous-$\infty$-catégorie des faisceaux pour la topologie $\cT$.
Dans le cas où $C=\Aff_k$ est la catégorie des $k$-schémas affines munie de la topologie étale sur un corps $k$, on appelle catégorie des \textsl{préchamps} (resp. des $\infty$-\textsl{champs}) la catégorie $\PreSh(\Aff_k)$ (resp. $\St_k=\Sh(\Aff_k)$).
L'inclusion $\Sh(C)\hra\PrSh(C)$ admet un adjoint $\PrSh(C)\ra \Sh(C)$, donnée par la faisceautisation \cite[6.2.2.7]{Lu1}.

On reprend la définition de \cite[1.2.4]{BKV}.
Soit (P) une classe de morphismes de $\infty$-champs $f:\cX\ra Y$ avec $Y\in\Aff_k$, stable par changement de base.

Un morphisme $f:\cX\ra\cY$ de $\infty$-champs est $(P)$-représentable, si pour tout changement de base $Y\ra \cY$, avec $Y$ un schéma affine, le changement de base $\cX\times_{\cY}Y\ra Y$ est dans (P).

On définit ainsi pour des morphismes de $\infty$-champs $f:\cX\ra\cY$ une notion de morphisme propre/de présentation finie/ind-fp-propre \footnote{i.e. ind-(propre de présentation finie)}/ind-représentable s'il est $(P)$-représentable, en prenant pour $(P)$, la classe des morphismes de schémas $\cX\ra Y$ qui sont propres/de présentation finie, la classe des morphismes $\cX\ra Y$ ind-fp-propres avec $\cX$ un ind-espace algébrique et la classe des morphismes $\cX\ra Y$ avec $\cX$ un ind-espace algébrique, pour $Y$ un schéma affine.

\bexa\label{exquot}
Pour un $\infty$-champ $\cX$ muni d'une action d'un ind-schéma en groupes $H$, on peut former l'$\infty$-champ quotient $[\cX/H]$.
\eexa
\blem\label{eqrep}
Soit une paire $(H,\cX)$ comme dans \ref{exquot}, alors :
\benumr
\item
Le flèche $\cX\ra [\cX/H]$ est ind-représentable.
\item
Si de plus $H$ est  ind-(quasi-affine), la flèche $\cX\ra [\cX/H]$ est ind-schématique.
\eenum
\elem
\brem
Par ind-représentable, on entend que pour tout morphisme $\Spec(A)\ra [\cX/H]$, le changement de base est représentable par un ind-espace algébrique.
\erem
\bpf
D'après \cite[1.2.6.(d)]{BKV}, $[\cX/H](A)$ classifie les paires $(E,\phi)$ où $\pi:E\ra\Spec(A)$ est un morphisme de $\infty$-champs et un $H$-torseur pour la topologie étale tel que $[E/H]\cong \Spec(A)$ et $\phi:E\ra \cX$ est un morphisme $H$-équivariant. 

On commence par voir que $E$ est représentable par un ind-espace algébrique.
Comme $\pi$ est un $H$-torseur pour la topologie étale, il existe $\Spec(A')\ra\Spec(A)$ étale, tel que $E$ se trivialise. En particulier, $E':=E\times_{\Spec(A)}\Spec(A')$ est représentable par un ind-schéma. D'après \cite[Lem. 3.12]{HR}, pour voir que $E$ est représentable par un ind-espace algébrique, il suffit de vérifier que $\pi:E\ra\Spec(A)$ a une diagonale schématique. 
Or $\pi$ est un $H$-torseur de telle sorte que $\Delta_{\pi}:E\ra E\times_{\Spec(A)}E\cong E\times_{k}H$ et la diagonale s'identifie à la section unité. Comme par changement de base, $E\times_{k}H\ra E$ est ind-schématique et ind-(quasi-séparé), d'après \cite[Tag. 01KT]{Sta}  la section unité est une immersion schématique quasi-compacte.

Il reste à montrer que $E$ est représentable par un ind-schéma. On l'écrit alors $E\simeq\colim E_{\al}$ où $E_{\al}$ est un espace algébrique qcqs. Or, comme $H$ est ind-(quasi-affine), $E'$ l'est aussi de telle sorte que $E'\simeq\colim E'_{\al}:= E_{\al}\times_{A}\Spec(A')$ et donc $E'_{\al}$ est représentable par un schéma quasi-affine. Par effectivité de la descente pour des morphismes quasi-affines \cite[Tag. 0247]{Sta}, $E_{\al}$ est quasi-affine  et donc $E$ est représentable par un ind-schéma, comme souhaité.
\epf

Pour les morphismes ind-fp-propres, on a le lemme suivant tiré de \cite[Rems. 1.2.9. (a),(b)]{BKV}:
\blem\label{loopeq}
\benumr
	\item 
	La propriété d'être ind-fp-propre est étale locale sur la base.
	\item
	Soit un morphisme d'$\infty$-champs $f:\cX\ra\cY$ ind-fp-propre et équivariant pour l'action d'un $\infty$-champ en groupes $H$, alors $f:[\cX/H]\ra [\cY/H]$ est ind-fp-propre.
\eenum
\elem
\bexa
Une illustration de cet énoncé général est la suivante; soit un $k$-ind-schéma $X$ et un groupe $G$ connexe réductif sur un corps $k$ tel que $\cL G$ agit sur $X$. Alors, on peut former les $\infty$-champs $[X/\clp G]$ et $[X/\cL G]$ et la flèche induite $[X/\clp G]\ra[X/\cL G]$ est ind-fp-propre, puisque la projection $X\times_k[\cL G/\clp G]\ra X$ est ind-fp-propre comme la grassmannienne affine $[\cL G/\clp G]$ l'est.
\eexa

Pour tout $\infty$-champ $\cX$, on a d'après \cite[Cor. 5.1.5.8]{Lu1}, $\cX\cong\colim\limits_{S\ra\cX} S$ avec $S$ affine; on associe alors son réduit $\cX_{red}\cong\colim\limits_{S\ra\cX}S_{red}$.
On a donc une flèche canonique $\cX_{red}\ra \cX$.

\'{E}tant donné un morphisme d'$\infty$-champs $f:\cX\ra\cY$, on dit que c'est une \textsl{équivalence topologique}, s'il induit une équivalence $f_{red}:\cX_{red}\stackrel{\sim}{\rightarrow}\cY_{red}$.
D'après \cite[1.5.1]{BKV}, cette notion est locale pour la topologie étale sur la base, stable par tirés-en-arrière, composition et passage au quotient.

Pour tout morphisme de $\infty$-champs $i:\cY\ra X$, on dit que $i$ est une \textsl{immersion topologiquement constructible} si pour tout schéma affine $X\ra \cX$, la composée :
\[(\cY\times_{\cX}X)_{red}\ra\cY\times_{\cX}X\ra X\]
est une immersion localement fermée de présentation finie.
\subsection{Fibration de Grothendieck-Springer affine}

\subsubsection{Grassmanniennes affines}
Soit un anneau commutatif $A$,  $G$ un schéma en groupes affine lisse sur $\Spec(A)$. On considère alors la grassmannienne affine de $G$, comme le quotient étale $\Gr_G=[\cL G/\clp G]$.
On a le lemme suivant tiré de \cite[3.10]{HR}:
\blem\label{rep}
Soit $G'\ra\Spec(A)$ un schéma en groupes affine lisse.
Alors pour tout sous-groupe fermé affine lisse $G\hra G'$ tel que le quotient fppf $G'/G$ est représentable par un $A$-schéma quasi-affine (resp. affine),  $\Gr_{G}\ra \Gr_{G'}$ est représentable par une immersion quasi-compacte (resp. immersion fermée).
De plus, si $G'=GL_n$, $\Gr_{GL_n}$ est représentable par un $A$-ind-schéma ind-projectif.
\elem
\brem\label{tho}
Il résulte de \ref{rep} et \cite[Thm. 9.4]{Alp} que si $G$ est réductif et admet un plongement dans $GL_n$, alors $\Gr_G$ est représentable par un ind-schéma ind-projectif. C'est par exemple le cas si $G$ est réductif déployé ou si $A$ est normal d'après \cite{Tho}.
\erem
De plus, si $G$ admet un Borel $B$, on pose $I=\ev^{-1}(B)$ et on introduit la variété de drapeaux affine $\Fl=\cL G/I$, également représentable par un ind-schéma ind-projectif. 
\medskip

Pour la suite, il est utile d'établir des lemmes de plongement d'un schéma en groupes dans un groupe linéaire avec quotient quasi-affine.

\blem\label{aff-quot}
Soit $B$ une $A$-algèbre de présentation finie, finie, libre. Soit $\pi:\Spec(B)\ra\Spec(A)$
\benumr
\remi
Soit $T$ un tore déployé sur $\Spec(B)$, alors la restriction à la Weil $J^2=\pi_*T$ admet un plongement linéaire avec un quotient quasi-affine.
\remi
On suppose de plus $\pi$ muni de l'action d'un groupe fini $\Gm$ d'ordre premier aux caractéristiques résiduelles.
Soit $G$ un schéma en groupes lisse sur $\Spec(B)$, alors $J^1=(\pi_*G)^{\Gm}$ est lisse et $J^2/J^1$ est quasi-affine.
\eenum
\elem
\bpf

(i) Quitte à considérer le produit de plongements et à faire le produit des espaces quotients, il suffit de traiter le cas $T=\bG_m$.
Soit $n$ le rang de $B$ comme $A$-module.
Dans ce cas, on a $J^2=\pi_{*}\bG_m=\Aut_{A-alg}(B)\subset GL_n=\Aut_{A}(B)$ et $J^2$ est un schéma en groupes lisse sur $\Spec(A)$. Montrons que $ GL_n/J^{2}$ est quasi-affine.
Comme $B$ admet une structure de $A$-algèbre, on a donc un morphisme $\xi: B\otimes_{A}B\ra\ B$. Comme par hyptohèse $B$ est un $A$-module libre de type fini, on peut voir $\xi$ comme un élément de $V=B^{\vee}\otimes_{A}B^{\vee}\otimes_{A}B$ sur lequel $ GL_n$ agit et donc $J^{2}$ est le stabilisateur de $\xi$. Le quotient $ GL_n/J^{2}$ est représentable d'après \cite[Exp.V, Thm. 10.1.2]{SGA3} et le quotient $ GL_n/J^{2}$ s'identifie à l'orbite de $\xi$ dans $V$. Ainsi, d'après \cite[II. 5. Prop.3.1]{DG}, elle est ouverte dans son adhérence, donc elle est quasi-affine.

(ii) Considérons l'immersion fermée $J^{1}\hra J^{2}=\pi_{*}G$ et montrons que le quotient $J^{2}/J^{1}$ est affine.
Tout d'abord, comme l'ordre de $\Gm$ est premier aux caractéristiques résiduelles, $J^{1}$ est un schéma en groupes lisse (\cite[sect. 15]{GKM}).
On fait donc agir $J^{2}$ par conjugaison sur $\prod\limits_{w\in \Gm}J^{2}\rtimes \Gm$ de telle sorte que $J^{1}$ est le stabilisateur du uplet $\xi=(1,w)_{w\in \Gm}$.
On applique alors \cite[Exp.V, Thm. 10.1.2]{SGA3} pour obtenir que $J^{2}/J^{1}$ est représentable. Comme de plus, il s'identifie à l'orbite $J^{2}.\xi$ dans $\prod\limits_{w\in \Gm}J^{2}\rtimes \Gm$, à nouveau d'après \cite[II. 5. Prop. 3.1]{DG}, elle est ouverte dans son adhérence, donc en particulier est quasi-affine.

\epf
\subsubsection{Une fibration en ind-schémas}
On rappelle les définitions de \cite[4.1]{BKV}. Soit un groupe connexe réductif déployé $G$ sur un corps $k$ algébriquement clos.
La flèche $\chi:\kg\ra\kc$ induit un morphisme au niveau des lacets, noté de la même manière $\chi:\cL\kg\ra\cL\kc$, on considère le ind-schéma des éléments compacts $\kC=\chi^{-1}(\clp\kc)$. C'est un sous-ind-schéma fermé de présentation finie de $\clp\kg$ et on considère l'ouvert des élements compacts génériquement réguliers semisimples:
\[\kC_{\bullet}=\chi^{-1}(\clp\kc-\clp\kD)\]
On prendra garde au fait que $\kC_{\bullet}$ n'est pas un ind-schéma au sens strict, car $\clp\kc-\clp\kD$ n'est pas quasi-compact. Néanmoins, si on l'écrit comme une union croissante $\clp\kc-\clp\kD=\bigcup_{d\in\NN}\clp\kc_{\leq d}$ où l'on borne la valuation du discriminant, alors les ouverts réciproques $\kC_{\leq d}$ sont bien des ind-schémas.

On forme alors le sous-ind-schéma fermé $\Fl\times\kC$:
\[\widetilde{\kC}=\{(g,\g)\in\Fl\times\kC,\ad(g)^{-1}(\g)\in\Lie(I)\}\]
La seconde projection donne une flèche:
\[f:\widetilde{\kC}\ra\kC,\]
que l'on appelle la fibration de Grothendieck-Springer affine. Elle est ind-fp-propre.

Pour tout $\g\in\kC_{\bullet}(k)$, la fibre:
\[X_{\g}=f^{-1}(\g)=\{g\in\Fl,\ad(g)^{-1}(\g)\in\Lie(I)\}\]
est appelée la fibre de Springer affine. C'est un ind-schéma tel que $X_{\g,red}$ est un $k$-schéma localement de type fini, de dimension finie et équidimensionnel d'après Kazhdan-Lusztig \cite[sect.2, Prop.1, sect. 4, Prop.1]{KL}.
Il est commode également d'introduire la variante sphérique:
\[\widetilde{\kC}_{K}=\{(g,a)\in\Gr\times\kC,\ad(g)^{-1}(\g)\in\clp\kg\}.\]
Ici, $K$ désigne le compact maximal $\clp G$. On a en fait un diagramme où tous les carrés sont cartésiens:
\begin{equation}
\xymatrix{\wkC\ar[r]^-{\inv}\ar[d]&[Lie(I)/I]\ar[d]\ar[r]^{\ev}&[\Lie(B)/B]\ar[d]\\\kC\ar[r]&[\clp\kg/\clp G]\ar[r]^-{ev}&[\kg/G]}
\label{D-Sprin}
\end{equation}
où les quotient sont faits pour la topologie étale, les actions sont les actions adjointes et les flèches verticales sont les projections canoniques.
On  introduit l'ouvert régulier $\widetilde{\kC}^{reg}=\inv^{-1}(\Lie(I)^{reg}/I)\subset \widetilde{\kC}$.
En général, au-dessus de $\kC$, $\widetilde{\kC}^{reg}$ peut avoir des fibres vides (e.g. pour $\g=0$). En revanche, si l'on se restreint au-dessus de $\kC_{\bullet}$, la flèche:
\[\widetilde{\kC}^{reg}\ra\kC_{\bullet}\]
est surjective d'après Kazhdan-Lusztig \cite[Cor. 1]{KL}.
De plus, si l'on considère la variante sphérique, alors pour tout $\g\in\kC_{\bullet}(k)$, la fibre de Springer affine sphérique $X_{\g,K,red}$ contient $X^{reg}_{\g,K,red}$ comme ouvert dense et le complémentaire est de dimension strictement inférieure d'après Ngô  \cite[Prop.3.10.1]{N}. 

\subsection{Symétries de la fibration}
Dans cette section, il s'agit de munir la fibration de Grothendieck-Springer affine $\wkC\ra\kC$ de l'action d'un ind-schéma en groupes, le Picard local. Cette construction est l'analogue local du champ de Picard qui agit sur la fibration de Hitchin \cite[sect. 4.3]{N}.
\subsubsection{Champ de Picard local}
Considérons l'espace d'arcs $\clp\kc$ de $\kc=\kt/W$, comme $\clp\kc=\Spec(A)$ est représentable par un schéma affine, on dispose d'une famille d'arcs universels:
\[e:D\hat{\times}\clp\kc=\Spec(A[[t]])\ra\kc,\]
et on peut considérer $\hat{J}=e^{*}J$ qui est un $A[[t]]$-schéma affine lisse, commutatif. On peut donc alors considérer sa grassmannienne affine:
\[\cP:=\Gr_{\hat{J}}\ra\clp\kc.\]
On l'appelle le champ de Picard local. En particulier, pour tout $a\in(\clp\kc)(k)$, on a $\cP_a(k)\cong J_a(k((t)))/J_a(k[[t]])$.
\blem\label{P-rep}
Le quotient $\cP$ est représentable par un ind-schéma de ind-présentation finie sur $\clp\kc$.
\elem
\bpf
On va appliquer à plusieurs reprises le lemme \ref{rep}. Tout d'abord, on commence par rappeler la description galoisienne du centralisateur régulier. D'après \cite[2.4.6-2.4.7]{N}, $J$ s'identifie canoniquement à un sous-schéma en groupes  ouvert de $J^{1}=\pi_*(T)^{W}$ et contient la composante neutre $J^{0}$ de $J^{1}$, de telle sorte que $J^{1}/J$ est un schéma en groupes fini sur $\kc$.
D'après \ref{rep}, il suffit donc de montrer la représentabilité de la grassmannienne affine de $J^1$ et cela se déduit alors du lemme \ref{aff-quot} et de \ref{rep} appliqué deux fois.
\epf
\subsubsection{Action sphérique du Picard local}
On a une flèche canonique $\chi:\kC\ra\clp\kc$ et on considère le changement de base $\cP_{\kC}=\chi^{*}\cP$.
On va construire une action de $\cP_{\kC}$ sur $\wkC$ et $\wkC_{K}$ au-dessus  de $\kC$.
\bprop\label{P-act}
On a une action de $\cP_{\kC}$ sur $\wkC$ au-dessus de $\kC$ qui commute à l'action de $\cL G$.
\eprop
\bpf
Soit une $k$-algèbre $A$ et $\g\in\kC(A)$ et $(g,\g)\in\wkC_K(A)$, avec $\ad(g)^{-1}\g\in\kg(A[[t]])$.
Posons $\g_0=\ad(g)^{-1}\g$ et $a=\chi(\g_0)$, d'après \eqref{carJ}, on a une flèche canonique au-dessus de $A((t))$:
\[\theta:J_{a}\ra C_{\g}\]
qui s'étend en un morphisme sur $\Spec(A[[t]])$:
\[\eta_{g}=\ad(g)^{-1}\circ\theta:J_{a}\ra C_{\g_0}\]
de telle sorte que :
\begin{equation}
\ad(g)^{-1}\circ\theta(J_{a}(A[[t]]))\subset C_{\g_0}(A[[t]])\subset G(A[[t]]).
\label{c-incl}
\end{equation}
Quitte à localiser pour la topologie étale, on peut supposer que $\cP(A)=J_{a}(A((t)))/J_{a}(A[[t]])$.
Pour tout $j\in J_{a}(A((t)))$, on pose alors:
 \begin{equation}
j.(g,\g)=(\eta_{g}(j).g,\g),
\label{act-eq}
\end{equation}
 qui agit sur $\wkC_K(A)$  et qui passe au quotient d'après \eqref{c-incl} en une action de $\cP_a(A)$ sur la fibre de Springer affine $\wkC_K(A)$.
Il reste à voir que cela commute à l'action de $\cL G$. L'action de $\cL G$ est donnée par:
\[k.(g,\g)=(kg,\ad(k)\g), k\in\cL G, (g,\g)\in\kC.\]
et cela vient du fait que $\eta_{kg}(j)=k\eta_{g}(j)k^{-1}$.
\epf

La section de Kostant induit un morphisme noté de la même manière $\eps:\clp\kc\ra \kC$ et on forme:
\[X_{\kc}:=\{(g,\g)\in\Fl\times\clp\kc, ad(g)^{-1}\eps(a)\in\Lie(I)\},\]
ainsi que sa variante sphérique $X_{K,\kc}$.

\bprop
Le morphisme $\wkC^{reg}_{K}\ra\kC$ est formellement principal homogène sous $\cP$, surjectif au-dessus de $\kC_{\bullet}$. De plus, si l'on tire $\wkC_{K}^{reg}$ par la section de Kostant $\eps:\clp\kc\ra \kC$, alors $\eps^{*}\wkC^{reg}_{K}\ra\clp\kc$ est le $\cP$-torseur trivial.
\eprop
\bpf
La surjectivité vient de \cite[Cor.1]{KL}
Montrons que l'on a un isomorphisme canonique:
\[\act:\cP\times_{\kC} \wkC_{K}^{reg}\ra \wkC_{K}^{reg}\times_{\kC}\wkC_{K}^{reg},\]
ce que nous entendons par le fait que la flèche est formellement principalement homogène sous $\cP$.
Soit un anneau $A$, montrons que tout élément $\g\in\kg^{reg}(A[[t]])$ est, localement pour la topologie étale sur $A$, conjugué à sa section de Kostant $\g_0=\eps(\chi(\g))$.
On considère alors le carré cartésien:
$$\xymatrix{E_{J}\ar[r]\ar[d]&G\times\kg^{reg}\ar[d]\\\Spec(A[[t]])\ar[r]^{(\g,\g_0)}&\kg^{reg}\times_{\kc}\kg^{reg}}$$
D'après \eqref{Jtors} et par changement de base, $E_{J}$ est un $J_a$-torseur sur $\Spec(A[[t]])$. D'après \cite[Thm. 2.1.7]{BC}, il est localement trivial pour la topologie étale sur $A$, donc localement sur $A$,  il existe un élément $k\in G(A[[t]])$ tel que $\ad(k)^{-1}\g=\g_0$.
On considère maintenant un triplet $(g_1,g_2,\g)\in \wkC_{K}^{reg}\times_{\kC}\wkC_{K}^{reg}(A)$.

Posons $\g_1=\ad(g_1)^{-1}\g$ et $\g_2=\ad(g_2)^{-1}\g$ dans $\kg^{reg}(A[[t]]])$ et $a=\chi(\g)$. D'après ce que l'on vient de voir, ils sont tous deux conjugués à $\eps(a)$. En particulier, $\eps(a)$ est $G(A((t)))$-conjugué à  $\g$ et $g_1$ et $g_2$ diffèrent par un élément de $J_{a}(A((t)))$, comme souhaité.
Maintenant si l'on tire par la section de Steinberg, $\eps^{*}\wkC^{reg}_{K}\ra\clp\kc$ admet une section donnée par $a\mapsto (1,\eps(a))$, donc c'est le $\cP$-torseur trivial.
\epf
\brem Il n'est pas clair que $\wkC^{reg}_{K}\ra\kC_{\bullet}$ soit un $\cP$-torseur, même sous la $h$-topologie, la raison à cela est bien que la flèche soit surjective et se factorise à travers une immersion ouverte quasi-compacte suivie d'un morphisme ind-fp-propre, la composée n'est pas nécessairement un $h$-recouvrement.
Une autre façon de le voir est qu'a priori ce morphisme n'est pas surjectif sur les $R$-points pour $R$ un anneau de valuation.
\erem
\subsubsection{Extension de l'action}\label{P-act2}
On veut maintenant relever l'action de $\cP$ à tout $\wkC$.
On commence par former le carré cartésien:
$$\xymatrix{\kC^{\sharp}\ar[d]\ar[rr]&&\kt\ar[d]\\\kC\ar[r]&\clp\kc\ar[r]&\kc}$$
Il résulte du diagramme commutatif 
$$\xymatrix{[\kb/B]\ar[d]\ar[r]&\kt\ar[d]\\[\kg/G]\ar[r]&\kc}$$
et de \eqref{D-Sprin} que la flèche $f:\wkC\ra\kC$ se factorise en:
\[f:\wkC\ra\kC^{\sharp}\ra\kC.\]
La proposition \ref{B-var} donne alors que $\cP\times_{\clp\kc}\kC^{\sharp}$ agit sur $\wkC$ au-dessus de $\kC^{\sharp}$.
De manière explicite, pour  un anneau $A$,  soit $(x,t_{B})\in \wkC(A)$ avec $x\in \wkC_K(A)$ d'image $a\in\kC(A)$ et $t_{b}\in[\kb/B](A)$, alors la flèche canonique au-dessus de $A[[t]]$ 
\[J_{a}\ra\Aut(x)\]
se factorise canoniquement modulo $t$ par :
\[J_{a}\ra Aut(t_b)\ra\Aut(x).\]
Et à nouveau elle commute à l'action de $\cL G$ sur $\wkC$.
\subsection{Description de la fibration au-dessus d'une stratification}
\subsubsection{Stratification par les valuations radicielles}
On a un diagramme commutatif:
\begin{equation}
\xymatrix{\kt\ar[dr]_{\pi}\ar[r]&\kg\ar[d]^{\chi}\\&\kc},
\label{Ccom}
\end{equation}
où la flèche horizontale est l'inclusion canonique.

Soient $m=\vert W\vert$, $\cO=k[[t]]$ et $F=k((t))$, $\co'=k[[t^{\frac{1}{m}}]]$, $F'=\Frac(\co')$, $\kt'=\Res_{\cO'/\co}(\kt\times_k\cO')$ et $\kc'=\Res_{\co'/\co}(\kc\times_k\cO')$. On rappelle que l'ordre de $W$ est premier à la caractéristique.
On définit alors $T_w=\Res_{\cO'/\co}(\kt\times_k\cO')^{w\mu_{m}}$ où $\mu_{m}$ agit sur $\cO'$ par l'action de Galois. D'après \cite[Lem. 15.3.1, 15.4.1]{GKM}, c'est un schéma en groupes lisse sur $\co$, qui est un tore sur $F$ et $\kt_w:=\Lie(T_w)=\Lie(\kt')^{w\mu_m}$.
D'après \cite[4.2, 4.3]{GKM2}, $T_{w}$ se déploie en $T$ sur $F'$ et on a un plongement $T_w\hra G$ unique à $G$-conjugaison près. De plus, tout tore maximal de $G$ est $G$-conjugué à un $T_w$, pour $w\in W$.
En particulier, tout élément $\g\in\kg^{rs}(F)$ est $G(F)$-conjugué à un élément de $\kt_{w}(F)$.
En considérant la composée $\kt_w\ra\kg\ra\kc$, on obtient une flèche canonique, indépendante des choix, définie sur $F$, $\pi_w:\kt_w\ra\kc.$
\blem\label{w-fin}
La flèche $\pi_w$ est finie.
\elem
\bpf
Par descente, il suffit de montrer que le morphisme est fini après extension de scalaires $F'/F$. Sur $F'$, $T_w$ se déploie en $T$ et la flèche $\pi_{w}\times_{F} F'$ s'identifie à la flèche $\pi:\kt\ra\kc$, qui est également finie, ce qu'on voulait.
\epf
Pour toute fonction $r:R\ra\frac{1}{m}\NN$ sur l'ensemble des racines, on peut considérer le sous-schéma localement fermé de présentation finie $\kt'_{r}$ où l'on fixe l'ensemble des valuations radicielles. En intersectant avec $\clp\kt_w$, on obtient ainsi la strate $\kt_{w,r}$.
De manière concrète, $\clp\kt_{w}$ classifie les séries $\sum\g_{i}t^{\frac{i}{m}}$ telles que $w^{-1}(\g_{i})=\xi^{i}\g_{i}$ et $\kt_{w,r}\subset\clp\kt_{w}$ est donnée par les équations:
\begin{equation}
d\al(\g_{i})=0~~ \text{si}~~0\leq i\sinf r(\al)~~ \text{et}~ d\al(\g_{i})\neq0~~ \text{si}~~ r(\al)=i.
\label{eqsplit}
\end{equation}
Ainsi, $\kt_{w,r}$ est donné par un nombre fini d'équations et d'inéquations linéaires, donc $\kt_{w,r}$ est fortement pro-lisse.
Pour $u\in W$ et une fonction $r$, on considère $u.r:R\ra\frac{1}{m}\NN$ donnée par $u.r(\al)=r(u^{-1}.\al)$. 
On a donc une action de $W$ sur les paires $(w,r)$  donnée par $u.(w,r)=(uwu^{-1},u(r))$.
On note $\kc_{w,r}\subset\clp\kc$ l'image de $\kt_{w,r}$. D'après \cite[Thm. 8.2.2]{GKM}, c'est un localement fermé de présentation finie, irréductible et fortement pro-lisse et  d'après \cite[3.3.4]{BKV}, on a une flèche finie étale $\kt_{w,r}\ra\kc_{w,r}$ de groupe $W_{w,r}$, où $W_{w,r}\subset W$ est le stabilisateur de $(w,r)$ pour l'action de $W$ sur les paires $(w,r)$.
Soit $\Lie(I)_{w,r}$ l'image inverse de $\kc_{w,r}$ par la flèche $\chi:\Lie(I)\ra\clp\kc$.
\blem\label{lci}
L'immersion $\Lie(I)_{w,r}\hra\Lie(I)$ est une immersion régulière de présentation finie.
\elem
\bpf
En effet, d'après \cite[Cor. 3.4.8]{BKV}, $\Lie(I)\ra\clp\kc$ est plat. De plus $\kc_{w,r}$ est lisse, donc $\kc_{w,r}\ra\clp\kc$ est une immersion régulière et d'après \cite[Tag. 067P]{Sta}, il en est de même de  $\Lie(I)_{w,r}\hra\Lie(I)$  par changement de base.
\epf

On définit $\kC_{w,r}=\kc_{w,r}\times_{\cL\kc^{rs}}\cL\kg^{rs}$ et son image inverse $\wkC_{w,r}\subset\wkC$.
On commence par un énoncé de structure sur le $\infty$-champ quotient $[\kC_{w,r}/\cL G]$.
Tout d'abord, d'après \cite[Lem. 4.1.8]{BKV}, on a un isomorphisme canonique:
\[\cL(G/T_{w})\times^{W_{w,r}}\kt_{w,\br}\cong\kC_{w,r}.\]
induit par $(g,t)\mapsto\ad(g).t$ et par composée on obtient donc un morphisme:
\begin{equation}
\phi:[\cL G/\cL T_{w}]\times^{W_{w,r}}\kt_{w,r}\ra \cL(G/T_{w})\times^{W_{w,r}}\kt_{w,r}.
\label{wrL}
\end{equation}
On a alors la proposition suivante:
\bthm\label{Tquot}
Pour tout $w\in W$, la flèche canonique de $\infty$-champs $\la:[\cL G/\cL T_w]\ra\cL(G/T_w)$ est un isomorphisme.
\ethm
\bpf
D'après \cite[Thm. 4.1.9]{BKV}, c'est un isomorphisme au niveau des réduits.
D'après \ref{loop}, comme $G/T_w$ est lisse, $\cL (G/T_w)$ est ind-placide et est ind-affine, on considère un sous-schéma fermé $\Spec(A)\ra\cL(G/T_w)$ avec $\Spec(A)$ placide.

On obtient alors par tiré-en-arrière un $T_w$-torseur $E_{A}\ra\Spec(A((t)))$. Il faut montrer qu'il est trivial localement pour la topologie étale sur $\Spec(A)$. Soit $\al\in H^{1}(A((t)),T_w)$ la classe associée.
Comme $\la_{red}$ est un isomorphisme, il existe un recouvrement étale $\Spec(B)\ra\Spec(A_{red})$ tel que $\al_{\vert B}\in H^{1}(B((t)),T)$ est triviale.

Comme $\Spec(A)$ est placide, $\Spec(A)_{red}\hra\Spec(A)$ est une immersion fermée de présentation finie d'après \cite[Cor. 1.4.5]{BKV}
et donc l'idéal $\nil(A)$ est nilpotent de type fini. Ainsi, d'après \cite[Exp. III, Cor. 6.8]{SGA1}, il existe un unique $\Spec(B')\ra\Spec(A)$ étale qui se réduit au-dessus de $A_{red}$ sur $B$.
Or la paire $(B'((t)),\nil(A).B'((t)))$ est hensélienne, de telle sorte que d'après \cite[Thm. 2.1.7]{BC}, on a $H^{1}(B'((t)),T_w)\cong H^{1}(B((t)),T_w)$. Ainsi, $\al_{\vert B'}=1$ et $\la$ est un isomorphisme.
\epf

\bcor\label{wrquot}
La flèche $\phi$  induit un isomorphisme:
\[[\kt_{w,r}/(W_{w,r}\rtimes\cL T_w)]\stackrel{\sim}{\rightarrow} [\kC_{w,r}/\cL G].\]
En particulier, on en déduit un isomorphisme $[\kt_{w,r}/(W_{w}\rtimes\cL (T_{w})_{red})]\stackrel{\sim}{\rightarrow} [\kC_{w,r}/\cL G]_{red}$.
\ecor
\bpf
D'après \ref{Tquot}, la flèche \eqref{wrL} est un isomorphisme, il suffit donc de passer au quotient par $\cL G$ pour obtenir l'isomorphisme souhaité.
Enfin, pour déduire l'identité au niveau des réduits, comme $\kt_{w,r}$ et $W_{w,r}$ sont réduits, $[\kt_{w,r}/(W_{w}\rtimes\cL T_{w})]_{red}\cong [\kt_{w,r}/(W_{w,r}\rtimes(\cL T_{w})_{red})]$.
\epf
On forme alors le carré cartésien:
\begin{equation}
\xymatrix{\ti{X}_{w,r}\ar[d]_{\ti{f}_{w,r}}\ar[r]&[\wkC_{w,r}/\cL G]=[Lie(I)_{w,r}/I]\ar[d]^{\ov{p}}\ar[r]&[\wkC/\cL G]=[\Lie(I)/I]\ar[d]\\\kt_{w,r}\ar[r]^-{\Upsilon_{w,r}}&[\kC_{w,r}/\cL G]\ar[r]&[\kC/\cL G]}.
\label{wr-diag}
\end{equation}
D'après \cite[Cor. 4.3.4]{BKV}, on a les assertions suivantes:
\benumr
\item
 $X_{w,r}=(\ti{X}_{w,r})_{red}\ra\kt_{w,r}$ est un morphisme de schémas, localement de présentation finie, muni de l'action d'un réseau $\La_{w}$.
\item 
On a $X_{w,r}\simeq\kt_{w,r}\times_{[\kC_{w,r}/\cL G]_{red}}[\wkC_{w,r}/\cL G]_{red}$ et $[\wkC_{w,r}/\cL G]_{red}\cong[\Lie(I)_{w,r,red}/I]$.
\item
Le quotient $[X_{w,r}/\La_{w}]$ est un espace algébrique, fp-propre sur $\kt_{w,r}$.
\eenum

\section{Rappels sur les résultats de \cite{BKV}}
Dans cette section, il s'agit de rappeler et de généraliser certains résultats établis par Bouthier-Kazhdan-Varshavsky dans \cite{BKV}. Les résultats principaux sont \ref{fond-spr2} ainsi que \ref{lusz} et \ref{P-act3} qui comparent l'action du groupe de Weyl affine construite dans \cite{BKV} avec celle de Lusztig et enrichit le faisceau de Grothendieck-Springer affine des symétries qui viennent du Picard local.
 \subsection{Catégorie de faisceaux}
\subsubsection{Définition}\label{l-ad}
Soit $\ell$ un premier différent de la caractéristique de $k$, $\AlgSp^{tf}_k$ (resp. $\AlgSp_k$) la catégorie des $k$-espaces algébriques de type fini (resp. qcqs).
\begin{itemize}
	\remi
Soit $\Catst$  l'$\infty$-catégorie des petites $\infty$-catégories, stables et $\bql$-linéaires avec des foncteurs exacts, i.e. qui préservent les colimites finies, comme morphismes. D'après \cite[1.1.4.4, 1.1.4.6]{Lu2}, elle contient toutes les petites colimites filtrantes et toute les petites limites.
\remi
Soit $\PrCat$ la sous-$\infty$-catégorie des  $\infty$-catégories présentables (\cite[5.5.0.1]{Lu1}), stables, $\bql$-linéaires, avec des $\infty$-foncteurs continus, i.e. les catégories sont stables par petites colimites et tous les foncteurs commutent aux petites colimites.
Elle contient toutes les petites limites et colimites d'après \cite[4.2.4.8, 5.5.3.13, 5.5.3.18]{Lu1}.
\end{itemize}
Tous les foncteurs qui apparaissent ici sont des $\infty$-foncteurs et les limites et colimites sont à prendre au sens homotopique.
Etant donné $\cC\in\Catst$, on peut former sa catégorie des ind-objets $\Ind(\cC)$ (\cite[5.3.5.1]{Lu1}) qui est stable (\cite[1.1.3.6]{Lu2}) et donc présentable (elle admet toutes les colimites finies et les colimites filtrantes, donc toutes les colimites). On a ainsi un foncteur naturel :
\[\Ind:\Catst\rightarrow\PrCat\]
donné par $\cC\mapsto \Ind(\cC)$, qui commute aux petites colimites filtrantes (\cite[5.3.5.10]{Lu1}, \cite[1.9.2]{DG} et \cite{Roz}).

Pour un $k$-espace algébrique de type fini $Y$, on dispose, d'après  Liu-Zheng (\cite{LZ1}, \cite{LZ2}), d'une $\infty$-catégorie $\cD_{c}(Y):=\cD^{b}_{c}(Y,\bql)$ dont la catégorie homotopique est $D^{b}_{c}(Y,\bql)$.
Dans la suite, on utilise librement le formalisme des six foncteurs pour les $\infty$-catégories (cf. \cite{LZ1},\cite{LZ2}).
Il y a un foncteur naturel:
\begin{equation}
\cD_{c}:(\AlgSp^{tf}_k)^{op}\rightarrow\Catst
\label{tfdef1}
\end{equation}
qui  à chaque morphisme $f:X\rightarrow Y$ associe $f^{!}:\cD_{c}(Y)\rightarrow\cD_{c}(X)$.
On définit alors $\cD:=\Ind\circ\cD_{c}$. En faisant l'extension de Kan à gauche et à droite, on obtient des foncteurs :
\begin{equation}
\cD_{c}:\PreSh(\AlgSp_k)\rightarrow\Catst
\label{kan1}
\end{equation}
\begin{equation}
\cD:\PreSh(\AlgSp_k)\rightarrow\PrCat
\label{kan2}
\end{equation}
On note $\cD_{\bullet}$ pour désigner les foncteurs $\cD$ et $\cD_c$.
De manière explicite si $X$ s'écrit comme une limite projective cofiltrante $X\simeq\varprojlim X_{\al}$ où les $X_{\al}$ sont des $k$-espaces algébriques de type fini avec des morphismes de transition affines, on a $\cD_{\bullet}(X)\simeq\colim_{f^!}\cD_c(X_{\al})$ et si $T=\colim T_{\al}$ est un ind-espace algébrique avec $T_{\al}\in\AlgSp_k$, on a $\cD_{\bullet}(T)\simeq\varprojlim_{f^!}\cD_{\bullet}(T_{\al})$.
D'après \cite[5.2.8]{BKV}, $\cD_{\bullet}$ est un faisceau pour la topologie étale, donc il se factorise par $\Sh(\AlgSp_k)$ via le morphisme de faisceautisation.
En particulier, comme le plongement $\iota:\Aff_k\hra\AlgSp_k$ induit une équivalence de $\infty$-catégories par tiré-en-arrière $\iota^{*}:\Sh(\Aff_k)\stackrel{\sim}{\rightarrow}\Sh(\AlgSp_k)$, on peut voir $\cD_{\bullet}$ comme un foncteur sur $\St_k$.

\subsubsection{Fonctorialité}\label{ss-fonct}
Pour tout morphisme de $\infty$-champs $f:\cX\ra\cY$, on dispose d'un foncteur $f^!$ pour $\cD_{\bullet}$.
De plus, pour toute équivalence topologique $f:\cX\ra\cY$, on a d'après \cite[5.3.6]{BKV} une équivalence de catégories:
\begin{equation}
f^{!}:\cD_{\bullet}(Y)\stackrel{\sim}{\rightarrow}\cD_{\bullet}(\cX).
\label{topeq}
\end{equation}
L'énoncé suivant est obtenu dans \cite[Prop. 5.3.7]{BKV}.
\bprop\label{Base}
Soit $f:\cX\ra\cY$ un morphisme ind-fp-propre d'$\infty$-champs.
Alors,  $f^!$ admet un adjoint à gauche $f_!$ et	pour tout carré cartésien entre $\infty$-champs:
$$\xymatrix{\ti{\cX}\ar[r]^{\ti{h}}\ar[d]_{\ti{f}}&\cX\ar[d]^{f}\\\ti{\cY}\ar[r]^{h}&\cY}$$
l'application de changement de base:
 \begin{equation}
\ti{f}_!\ti{h}^{!}\ra h^{!}f_!
\label{basemor}
\end{equation} est une équivalence.
\eprop

Pour la suite, on a besoin de considérer une classe intermédiaire entre les ind-schémas et les $\infty$-champs, qui contient les objets tels que $[\cL X/\cL G]$ pour $G$ connexe réductif sur $k$ qui agit sur $X\in\Aff_k^{tf}$.
Tout d'abord, pour toute immersion topologiquement constructible de $\infty$-champs $\eta:\cY\hra\cX$, on dispose d'un foncteur $\eta_*$ (\cite[5.4.4]{BKV}).
On s'intéresse à l'existence d'un adjoint à gauche.
\bdefi\label{glu}
Soit $\cX$ un $\infty$-champ, on dit qu'il satisfait le recollement si pour toute immersion topologiquement constructible $\cY\hra\cX$, le foncteur $\eta_{*}$ admet un adjoint à gauche $\eta^*$.
\edefi 
\brems\remi
\label{glustk}
Les $\infty$-champs suivants satisfont le recollement:
\benumr
\item
Si $\cX$ satisfait le recollement alors pour toute immersion topologiquement constructible $\eta:\cY\hra\cX$, $\cY$ satisfait le recollement (\cite[Lem. 5.5.5.(a)]{BKV}).
\item
Si $\cX\simeq\colim \cX_{\al}$ où chaque $\cX_{\al}$ admet le recollement et les morphismes de transition sont des immersions ouvertes quasi-compactes (resp. des immersions fermées de présentation finie) alors $\cX$ satisfait le recollement (\cite[Lem. 5.5.5.(b)]{BKV}).
\item
D'après \cite[5.5.6]{BKV}, on obtient que tout espace algébrique placide satisfait le recollement ainsi que tout ind-espace algébrique ind-placide d'après (ii).
\item
Les quotients $[X/H]$ d'un ind-espace algébrique ind-placide $X$ par un groupe ind-placide $H$, i.e. un objet en groupe dans la catégorie des ind-espaces algébrique ind-placides satisfont le recollement, d'après \cite[Prop. 5.5.7]{BKV}. Cela s'applique donc à des quotients  $[\cL Y/\cL G]$ pour $G$ affine lisse sur $k$ qui agit sur $Y\in\Aff_k^{tf}$.
\remi
Pour toute immersion topologiquement constructible $\eta:\cY\hra\cX$ de $\infty$-champs satisfaisant le recollement, $\nu_*$ préserve les constructibles \cite[Lem. 5.4.1, 5.4.8.(c)]{BKV}.
\eenum
\erems
\subsection{t-structure}
\subsubsection{Stratification constructible}
Soit un $\infty$-champ $\cX$, on considère une collection $(\cX_{\al})_{\al\in\cI}$ de sous-champs topologiquement constructibles tels que $\cX_{\al}\cap\cX_{\beta}=\emptyset$ pour tout $\al\neq\beta$.
Pour tout sous-champ $\cX'\subset\cX$, on pose $\cI_{\cX'}=\{\al\in\cI, \cX_{\al}\subset\cX'\}$.
On dit que $\cX'$ est $\cI$-\textsl{adapté} si pour tout $\al\in\cI\backslash \cI_{\cX'}$, $\cX\cap\cX_{\al}=\emptyset$.

\bdefi
Avec les notations de ci-dessus, on dit que  les $\{\cX_{\al}\}$ forment:
\begin{enumerate}
	\item 
	une stratification \textit{finie constructible} s'il existe une suite finie croissante d'ouverts quasi-compacts $(\cU_{i})_{1\leq i\leq n}$ telle que l'on a un plongement $\cX_{\al_i}\hra\cU_{i}-\cU_{i-1}$ qui est une équivalence topologique.
	\item
		une stratification \textit{bornée constructible} si $\cX$ s'écrit comme une colimite filtrante de sous-champs ouverts quasi-compacts $\cX\cong\colim_{U\in\J}\cX_{U}$ où chaque $\cX_{U}$ est $\cI$-adapté et les $\{\cX_{\al}\}_{\al\in\cI_{\cX_{U}}}$ forment une stratification finie constructible.
\end{enumerate}
\edefi
Dans la suite, on dit qu'un $\infty$-champ est \textsl{placidement stratifié} s'il admet une stratification bornée constructible par des $\cX_{\al}$ qui sont des champs placides.

\subsubsection{Le cas des champs placidement stratifiés}\label{t-codim}
On a besoin d'introduire une classe assez vaste d'objets pour lesquels on peut parler de $t$-structure et adaptée à la présentation de $\cD$ à l'aide de $f^!$. Pour cette $t$-structure, le faisceau dualisant sur un schéma lisse sera pervers auto-dual.
\medskip

Pour un $k$-schéma de type fini équidimensionnel $X$, on dispose de la $t$-structure !-adaptée où l'on déclare qu'un objet $K\in\cD_c(X)$ est pervers si et seulement si $K[-\dim(X)]$ est pervers pour la $t$-structure perverse standard.
Cette $t$-structure se comporte de la façon suivante par rapport à la dualité de Verdier:
\begin{equation}
 K\in \mathstrut^{p}\cD^{\geq 0}(X)\Leftrightarrow\bD(K)\in \mathstrut^{p}\cD^{\leq 2\dim(X)}(X).
\label{dulcar}
\end{equation}

Pour un $k$-schéma de type fini général, il admet une stratification canonique par des schémas équidimensionnels et on obtient la $t$-structure sur $X$ en recollant les $t$-structures sur chacune des strates équidimensionnelles (\cite[6.2.3]{BKV}).

On considère ensuite un $k$-schéma affine placide $X$; d'après \cite[Prop. 6.3.1]{BKV}, il existe une unique t-structure sur $\cD_c(X)$ telle que pour tout morphisme fortement pro-lisse $f:X\ra Y$ avec $Y$ est de type fini, $f^!$ est $t$-exact où l'on munit $\cD_c(Y)$ de la $t$-structure $!$-adaptée.
Comme pour un schéma placide affine, on a $\cD=\Ind(\cD_c)$, la $t$-structure s'étend naturellement à $\cD$ (cf. \cite[6.3.2.(a)]{BKV}).

Notant $\cD_{\bullet}$ pour désigner indifférement $\cD$ ou $\cD_c$, si $\cX$ est un champ placide, en considérant un atlas placide $X$, on montre qu'il existe une unique $t$-structure sur $\cD_{\bullet}(\cX)$ telle que pour tout morphisme lisse $f:X\ra\cX$ avec $X$ affine placide, $f^!$ est $t$-exact  (\cite[sect. 6.3]{BKV}).

Enfin la dernière étape consiste à considérer des $\infty$-champs placidement stratifiés. Soit $(\cX,(\cX_{\al})_{\al\in\cI})$ un champ placidement stratifié qui satisfait le recollement des faisceaux.
On se donne une perversité $p_\nu$ adaptée à la stratitification, i.e. une fonction $p:\cX\ra\bZ$ constante sur chaque strate, soit une collection d'entiers $(\nu_{\al})_{\al\in\cI}$. On suppose  pour simplifier la stratification bornée.
D'après ci-dessous pour tout $\al\in\cI$, on dispose d'une $t$-structure naturelle sur chaque $\cD(\cX_{\al})$ et comme $\cX$ satisfait le recollement on a des foncteurs $\eta_{\al}^{*},\eta_{\al}^!:\cD(\cX)\ra\cD(\cX_{\al})$. Il est important ici de considérer $\cD$ plutôt que $\cD_c$ pour disposer de tels foncteurs.
D'après \cite[Prop. 6.4.2]{BKV}, subordonnée à $p_\nu$ et si la stratification est \textsl{bornée}, il existe alors une unique $t$-structure $(\mathstrut^{p}\cD^{\leq 0}(\cX),\mathstrut^{p}\cD^{\geq 0}(\cX))$ sur $\cD(\cX)$ telle que:
\begin{equation}
\cD^{\geq 0}(\cX)=\{K\in\cD(\cX),\forall~\al\in\cI, \eta_{\al}^{!}K\in\mathstrut^{p}\cD^{\geq -\nu_{\al}}(\cX_{\al})\},
\label{t1}
\end{equation}
\begin{equation}
\cD^{\leq 0}(\cX)=\{K\in\cD(\cX),\forall~\al\in\cI, \eta_{\al}^{*}K\in\mathstrut^{p}\cD^{\leq -\nu_{\al}}(\cX_{\al})\}.
\label{t2}
\end{equation}
Si la stratification n'est pas finie, on a également une $t$-structure avec seulement \eqref{t1} qui est vérifiée.
Dans les cas qui nous intéressent la perversité $p_\nu$ est essentiellement celle donnée par \og la codimension\fg~ des strates $\cX_{\al}$ dans $\cX$.

Enfin, en conservant les hypothèses de ci-dessus, si $j:\cU\hra\cX$ est un ouvert quasi-compact, d'après \cite[Lem. 6.4.5]{BKV}, $\cU$ est placidement stratifié pour la stratification induite et satisfait le recollement de telle sorte que $j^!$ est $t$-exact. De plus, pour tout $K\in\cD(\cU)$, on peut définir (cf. \cite[Cor. 6.4.10, Lem. 6.4.11]{BKV}) un prolongement intermédiaire $j_{!*}K$  avec $j^{*}j_{!*}K=K$ et tel que :
\benumr
	\item 
		Si la stratification est finie et que $j$ est adapté à la stratification, $j_{!*}K$ est l'unique extension perverse $\ti{K}$ de $K$ tel que pour tout $\al\in\cI\backslash\cI_{\cU}$, $\nu_{\al}^{*}\ti{K}\in\cD^{\leq -\nu_{\al}-1}(\cX_{\al})$ et $\nu_{\al}^{!}\ti{K}\in\cD^{\geq -\nu_{\al}+1}(\cX_{\al})$.
	\item
	$\End(K)=\End(j_{!*}K)$.
\eenum
\subsubsection{Systèmes locaux renormalisés}
Soit $X$ un $k$-schéma de type fini lisse, soit $\Loc_{\bql}(X)$ la catégorie des systèmes locaux  sur $X$.
On note $\Loc^{ren}(X)=\{\cL[2\dim_X],\cL\in\Loc_{\bql}(X)\}$ où $\dim_X$ est la fonction localement constante de dimension, de telle sorte que comme $X$ est lisse, on a l'inclusion:
\[\Loc_{\bql}^{ren}(X)\subset\Perv^{!}(X)\]
 où $\Perv^{!}(X)$ désigne la catégorie des faisceaux pervers pour la $!$-t-structure.
En particulier, pour tout morphisme lisse $f:Y\ra X$ de $k$-schémas lisses, on a donc un foncteur 
\[f^{!}:\Loc^{ren}(X)\ra\Loc^{ren}(Y).\]
Soit un $k$-schéma placidement présenté $X$ lisse, on peut donc former :
\[\Loc^{ren}(X)=\colim_{X\ra Y}\Loc^{ren}(Y),\]
où $X\ra Y$ parcourt les morphismes fortement pro-lisses avec $Y$ lisse.
Par construction, $\Loc^{ren}(X)$ est fonctorielle par rapport aux morphismes fortement pro-lisses entre $k$-schémas placidement présentés lisses et est locale pour la topologie étale.
De la sorte, on obtient donc pour tout champ placide $\cX$ lisse une catégorie $\Loc^{ren}(\cX)$
et l'on a par $t$-exactitude de $f^!$ (\cite[Prop. 6.3.1-6.3.3]{BKV}) pour des morphismes fortement pro-lisses, une inclusion:
\[\Loc^{ren}(\cX)\subset\Perv^{!}(\cX).\]

\blem\label{t-lci}
Soit $i:\cH\ra\cX$ une immersion localement fermée de présentation finie  entre champs placides de codimension $d$, avec $\cX$ lisse. Alors, pour tout $K\in\Loc^{ren}(\cX)$, on a $i^{*}K\in\mathstrut^{p}\cD^{\leq-2d}(\cH)$ et $i^{!}K\in\mathstrut^{p}\cD^{\geq 0}(\cH)$.
\elem

\bpf
Le même dévissage que \cite[Preuve 6.3.5.(c)]{BKV} nous ramène immédiatement au cas où $\cX=X$ est un $k$-schéma lisse de type fini et $i:H\hra X$ une immersion  de codimension $d$ entre $k$-schémas de type fini. Quitte à passer aux composantes connexes, on peut supposer $X$ connexe.
Dans ce cas, on a $K=\cL[2\dim X]$ pour un système local $\cL$ sur $X$, d'où $i^{*}K= i^{*}L[2\dim X]\in\mathstrut^{p}\cD^{\leq-2d}(\cH)$.

Pour l'assertion duale, en vertu de \eqref{dulcar}, il suffit de vérifier que $\bD(i^{!}\cL[2\dim X])\in\mathstrut^{p}\cD^{\leq 2\dim(H)}(H)$.
Or, comme $X$ est lisse, on a $\bD(i^{!}\cL[2\dim X])=i^{*}\cL^{\vee}$ et l'assertion est claire.
\epf
\subsubsection{Application à la fibration de Grothendieck-Springer affine}
On veut appliquer la théorie générale à la fibration:
\[f:\wkC\ra\kC.\]
La $t$-structure n'existe qu'une fois que l'on a divisé par l'action adjointe de $\cL G$. 
Notons de la même manière la flèche au niveau des quotients:
\[f:[\wkC/\cL G]\ra[\kC/\cL G].\]
Dans ce cas, on a alors $[\wkC/\cL G]\cong [\Lie(I)/I]$ qui est déjà un champ placide lisse. On dispose donc d'une $t$-structure sur $\cD([\Lie(I)/I])$ pour laquelle $\omega_{\Lie(I)/I}$ est pervers.
La situation est plus délicate pour le quotient $[\kC/\cL G]$ et on doit dans un premier temps se restreindre à l'ouvert des génériquement réguliers semisimples $\kC_{\bullet}$. Plus précisément, pour tout $m\in\NN$, on considère l'ouvert $\kC_{\leq m}$  de telle sorte que $[\kC_{\leq m}/\cL G]$ est un $\infty$-champ  qui satisfait le recollement avec comme stratification finie placide les $[\kC_{w,r}/\cL G]_{red}$ d'après \cite[Thm. 4.1.12, Lem. 7.1.2]{BKV}.
On commence par donner un définition de la notion de codimension pour des sous-champs constructibles de $[\kC_{\bullet}/\cL G]$.
Tout d'abord, d'après \cite[Lemmes 2.2.4, 2.2.5]{BKV} pour tout schéma placide, on a une notion de codimension. 
Dans le cas plus simple où $X\simeq\varprojlim X_{i}$ est placidement présenté avec des flèches de transition surjectives lisses (ce qui est le cas de $\clp Y$ si $Y$ est lisse de type fini) alors tout fermé de présentation finie $Z\subset X$ se descend à un cran fini en $Z_i\subset X_i$ et alors on pose $\codim_{X}(Z)=\codim_{X_i}(Z_i)$ et l'on vérifie que cela ne dépend d'aucun choix.
On veut maintenant une notion de codimension pour les localement fermés de présentation finie de $[\kC_{\bullet}/\cL G]$.
\bdefi
Pour tout $m\in\NN$ et tout localement fermé de présentation  finie $Z\hra \kC_{\leq m}$ $\cL G$-équivariant, on définit :
\[\codim_{\kC_{\leq m}}(Z)=\delta_{Z}+\codim_{\clp\kg}(Z\cap\clp\kg)\]
où $\delta_{Z}$ est le plus petit entier tel que $\kC_{\delta}\cap Z \neq\emptyset$.
\edefi
En particulier, d'après \cite{GKM}, pour toute paire $(w,r)$, on a la formule suivante de codimension:
\begin{equation}
\nu_{w,r}:=\codim_{\kC_{\bullet}}(\kC_{w,r})=2\delta_{w,r}+c_{w,r}+d(w,r).
\label{wrcodim}
\end{equation}
avec $c_{w,r}=\dim_{k}\kt-\dim_{k}\kt^{w}$, $\delta_{w,r}=\frac{d_{w,r}-c_{w,r}}{2}$ pour $d_{w,r}=\sum\limits_{\al\in R}r(\al)$ et $d(w,r)$ la codimension de $\kt_{w,r}$ dans $\kt_{w}$.
En particulier, pour toute paire $(w,r)$ telle que $\kC_{w,r}\subset\kC-\kC_{\leq 0}$, on a (\cite[4.4.2.(c)]{BKV}):
\begin{equation}
2\delta_{w,r}\sinf \nu_{w,r}.
\label{wrcodim2}
\end{equation} 

Soit $j_{rs}:[\kC_{rs}=\kC_{\leq 0}/\cL G]\hra[\kC_{\bullet}/\cL G]$, le théorème principal établi par Bouthier-Kazhdan-Varshavsky est le suivant \cite[Thm. 7.1.4]{BKV}:
\bthm\label{fond-spr}
Notons $\cS_{\bullet}=f_!\omega_{[\wkC_{\bullet}/\cL G]}\in\cD([\kC_{\bullet}/\cL G])$.
On équipe $\cD([\kC_{\bullet}/\cL G])$ de la perversité donnée par la codimension, alors $\cS_{\bullet}$ est pervers, on a
$\cS\cong j_{rs,!*}j_{rs}^{*}\cS$. De plus, on a un isomorphisme canonique d'algèbres $\End(\cS_{\bullet})\cong\bql[\widetilde{W}]$.
En particulier, $\cS_{\bullet}$ est canoniquement muni d'une action de $\widetilde{W}$.
\ethm
En fait, on peut montrer un énoncé plus fort qui va inclure d'autres \og systèmes locaux\fg~ que le dualisant.
Pour toute paire $(w,r)$, soit $\ti{i}_{w,r}:[\Lie(I)_{w,r}/I]\ra[\Lie(I)_{\bullet}/I]$ l'inclusion et $b_{w,r}=\codim_{\Lie(I)}(\Lie(I)_{w,r})$. Il résulte de \cite[Prop. 3.4.2]{BKV} et de \eqref{wrcodim}, que l'on a :
\begin{equation}
b_{w,r}=\delta_{w,r}+c_{w,r}+d(w,r), \nu_{w,r}\leq 2b_{w,r}.
\label{wrbis}
\end{equation}
\bthm\label{fond-spr2}
On considère $f:[\wkC_{\bullet}/\cL G]\cong[\Lie(I)_{\bullet}/I]\ra[\kC_{\bullet}/\cL G]$ avec les perversités de \ref{fond-spr}. Alors pour tout faisceau $K\in\Loc^{ren}([\Lie(I)_{\bullet}/I])$, $f_{!}K$ est pervers et on a $f_{!}K\cong j_{rs,!*}j^{*}(f_!K)$.
\ethm
\brems
\remi
L'énoncé vaut plus généralement pour un morphisme ind-fp-propre petit (ou semi-petit) au sens de \cite[2.4.9]{BKV}, mais les conditions étant techniques à formuler, on se contente de cette formulation pour la fibration de Grothendieck-Springer affine, qui est le principal exemple.
\remi
Dans l'optique de construire des caractères de profondeur nulle, il est nécessaire de pouvoir considérer des systèmes locaux plus généraux que seulement le dualisant, ainsi que cela apparaît dans \cite[sect. 3.5]{BeKV}.
\erems
\bpf
(i) L'argument suit \cite[Thm. 6.5.3]{BKV}. Posons $\cX=[\Lie(I)_{\bullet}/I]$ et $\cY=[\kC_{\bullet}/\cL G]$ et notons $(\cY_{\al})$ (resp. $(\cX_{\al}))$ les stratifications induites pas les $\kC_{w,r}$.
Soit $K\in\Loc^{ren}([\Lie(I)_{\bullet}/I])$, notons $\cS_{K}=f_{!}K$. On commence par montrer que $\cS_{K}\in\mathstrut^{p}\cD^{\geq 0}$.
En vertu de la description \eqref{t1}, il suffit de  montrer que pour toute strate $(w,r)$, on a:
\begin{equation}
 i_{\al}^{!}\cS_{K}\in\mathstrut^{p}\cD^{\geq-\nu_{\al}}(\cY_{\al}) .
\label{p-cond}
\end{equation}
avec une inégalité stricte  pour les strates $\kC_{w,r}$ contenues dans $\kC_{\bullet}-\kC_{rs}$ et $i_{\al}:\cY_{\al}\ra\cY$.
Par définition, on a un carré cartésien:
$$\xymatrix{\cX_{\al}\ar[d]_{f_{\al}}\ar[r]^{\ti{i}_{\al}}&\cX\ar[d]^{f}\\\cY_{\al}\ar[r]^{i_{\al}}&\cY}.$$
 Par changement de base \cite[5.3.7, 5.5.4]{BKV}, on a :
\[i_{\al}^{!}\cS_{K}\cong f_{\al,!}\ti{i}_{\al}^{!}K.\]
D'après \ref{lci}, $\ti{i}_{w,r}$ est une immersion régulière, donc de pure codimension donc d'après \ref{t-lci}, $\ti{i}_{\al}^{!}K\in\mathstrut^{p}\cD^{\geq 0}$  et $f_{\al,!}[-2\delta_{\al}]$ est $t$-exact à gauche par \cite[Lem. 6.3.6]{BKV}, d'où $f_{\al,!}\ti{i}_{\al}^!L\in\cD^{\geq-2\delta_{\al}}(\cY_{\al})\subset \cD^{\geq-\nu_{\al}}(\cY_{\al})$, où la dernière inégalité vient de \eqref{wrcodim} avec inégalité stricte d'après \eqref{wrcodim2} pour les strates contenues dans $\kC_{\bullet}-\kC_{rs}$.

Montrons maintenant que $\cS_{K}\in\mathstrut^{p}\cD^{\leq 0}$, ainsi que l'assertion sur l'extension intermédiaire.
Par changement de base \cite[Prop. 5.3.7]{BKV}, $i_{\al}^{*}\cS_{K}\cong f_{\al,!}\ti{i}_{\al}^{*}K$ et d'après \cite[6.3.5. (b)]{BKV}, $f_{\al}^{!}$ est $t$-exact à gauche donc par adjonction $f_{\al,!}$ est $t$-exact à droite. En particulier, en utilisant \ref{t-lci} et \eqref{wrbis}, on obtient que $i_{\al}^{*}\cS_{K}\in\mathstrut^{p}\cD^{\leq-2b_{\al}}(\cY_{\al})\subset\cD^{\leq-\nu_{\al}}(\cY_{\al})$, avec inégalité stricte d'après \eqref{wrcodim2} pour les strates contenues dans $\kC_{\bullet}-\kC_{rs}$, comme voulu.
\epf
\subsection{Action de Lusztig}
\subsubsection{Combinatoire affine}
Dans \cite[5.4]{Lu}, Lusztig construit une action de $\widetilde{W}$ sur l'homologie des fibres de Springer affines. A priori, elle diffère de celle construite par le théorème \ref{W-constr}. L'objet de cette section est de comparer les deux.

Soit $R^{\vee}$ l'ensemble des coracines, on note $W_{aff}=\bZ R^{\vee}\ltimes W$ où $\bZ R^{\vee}\subset X_*(T)$ le réseau des coracines. On note $\Delta_{aff}$ l'ensemble des racines affines. Le groupe $W_{aff}$ est un groupe de Coxeter engendré par $\Delta_{aff}$, i.e. c'est le groupe engendré par les $s_{\al}$ avec $\al\in\Delta_{aff}$ et les relations $(s_{\al}s_{\beta})^{m_{\al\beta}}=1$.
On a une suite exacte:
\[1\ra W_{aff}\ra\widetilde{W}\ra\Omega\ra 1,\]
où $\Omega$ est fini abélien. Cette suite est scindée et $\Omega$ consiste en le stabilisateur de $\Delta_{aff}$ pour l'action de $\widetilde{W}$ sur l'ensemble des racines affines. En particulier, on dispose d'une action de $\Omega$ sur $\Delta_{aff}$. 
 Ainsi, $\widetilde{W}$ est engendré par $W_{aff}$ et $\Omega$ avec en plus des relations de Coxeter, les relations $usu^{-1}=u(s)$ pour $u\in\Omega$ et $s\in W_{aff}$.
On a également que $\Omega$ s'identifie à $N_{\cL G}(I)/I$.
 On appelle parahorique standard  tout sous-groupe fermé de présentation finie $I\subset P\subset\cL G$.  Soit $P^{+}$ son radical pro-unipotent, $L_{P}=P/P^{+}$ le quotient réductif maximal et $W_{P}\subset W_{aff}$ le groupe de Weyl associé.
Les parahoriques standards sont en bijection avec les sous-ensembles stricts de $\Delta_{aff}$ via la flèche $P\mapsto\Delta_{P}$ où $\Delta_{P}\subset\Delta$ est le sous-ensemble des racines simples associées à $L_P$.
Enfin, soit $B_{P}=I/P^{+}$, c'est un Borel de $L_{P}$.
\subsubsection{Une autre action de $\widetilde{W}$}

On commence par construire une action de $\widetilde{W}$ sur $\cS$ au niveau de la catégorie homotopique. On a un carré cartésien:
$$\xymatrix{[\Lie(I)/I]\ar[d]_{\pi_{P}}\ar[r]^-{\ev}&[\Lie(B_P)/B_P]\ar[d]^{\pi_{P}}\\[Lie(P)/P]\ar[r]^-{\ev}&[\Lie(L_P)/L_{P}]}$$
où la flèche $\pi_{P}$ est projective et petite par la théorie de Springer usuelle \cite[Cor.10.4]{KW}.
Par changement de base propre, on a donc:
\[\pi_{P,!}\omega_{[\Lie(I)/I]}\cong\ev^{!}\pi_{P,!}\omega_{[Lie(B_{P})/B_{P}]}.\]
Par la théorie de Springer usuelle \cite[Cor.10.5]{KW},$\pi_{P,!}\omega_{[Lie(B_{P})/B_{P}]}$ est muni d'une action de $W_{P}$ et donc également $\cS=f_{!}\omega_{[\Lie(I)/I]}=f_{P,!}(\pi_{P,!}\omega_{[\Lie(I)/I]})$ avec $f_{P}:[\Lie(P)/P]\ra[\kC/\cL G]$.
De plus, on vérifie que pour toute paire $Q\subset P$ l'action de $ W_{Q}$ sur $\cS$ s'obtient comme la restriction de $W_{P}$ sur $\cS$. En particulier, elle engendre une action de  $W_{aff}$ sur $\cS$.
Il s'agit maintenant de l'étendre en une action de $\widetilde{W}$. En utilisant l'identification $\Omega\cong N_{\cL G}(I)/I$, on obtient une action canonique par multiplication à droite de $\Omega$ sur $\Fl$. On note $\omega\mapsto R_{\omega}$ cette action.
Comme $\ad(g)I=\ad(g\omega) I$  pour tout $\omega\in\Omega$, on a aussi une action de $\Omega$ sur $\wkC$ au-dessus de $\kC$ qui commute à l'action de $\cL G$, d'où une action de $\Omega$ à gauche sur $\cS$:
\[R_{\omega,*}^{-1}:\cS\ra\cS.\]
De plus, l'action de $\omega\in\Omega$ envoie un parahorique $P$ sur $\omega^{-1} P\omega$ et on a un carré commutatif:
\begin{equation}
\xymatrix{[\Lie(I)/I]\ar[d]_{\pi_{P}}\ar[r]^{R_{\omega}}&[\Lie(I)/I]\ar[d]^{\pi_{\omega^{-1} P\omega}}\\[\Lie(P)/P]\ar[r]^-{R_{\omega}}&[\Lie(\omega^{-1} P\omega)/\omega^{-1} P\omega]}
\label{omega}
\end{equation}
En particulier, l'action de $\Omega$ entrelace l'action de $W_{P}$ et $W_{\omega^{-1} P\omega}$ sur $\cS$ via l'isomorphisme
$\ad(\omega):W_{\omega^{-1} P\omega}\stackrel{\sim}{\rightarrow}W_{P}$.
En utilisant la présentation par générateurs et relations de $\widetilde{W}$, on obtient une action de 
$\widetilde{W}$ sur $\cS$ dans la catégorie homotopique de $Ho(\cD([\kC/\cL G]))$ de $\cD([\kC/\cL G])$.
Il s'agit maintenant de la relever à l'$\infty$-catégorie, on se restreint alors à la sous-catégorie $\cD([\kC_{\bullet}/\cL G])$, dans laquelle, d'après \ref{fond-spr}, $\cS_{\bullet}$ est pervers, donc dans le coeur d'une $t$-structure, on a en particulier d'après \cite[Rmq. 1.2.1.12]{Lu2}, une équivalence $\Perv([\kC_{\bullet}/\cL G])\cong Ho(\Perv([\kC_{\bullet}/\cL G]))$ et l'action se relève.
On appelle cette action, l'action de Lusztig.

\bthm\label{lusz}
L'action de Lusztig sur $\cS_{\bullet}$  est la même que celle construite par le théorème \ref{fond-spr}.
\ethm
\bpf
Comme $\cS_{\bullet}$ est pervers, il suffit de vérifier l'énoncé dans la catégorie homotopique. De plus, d'après \ref{fond-spr}, on a $\End(\cS_{\bullet})=\End(\cS_{rs})$, de telle sorte qu'il suffit de vérifier l'égalité au-desssus de $\kC_{rs}$. Dans ce cas, on va voir que les deux actions proviennent d'une action de $\widetilde{W}$ sur la fibration.
D'après \cite[4.2.4]{BKV} la fibration $f_{rs}:\wkC_{rs}\ra\kC_{rs}$ est, à une équivalence topologique près, un $\widetilde{W}$-torseur qui fournit l'action de \ref{fond-spr}. De plus, pour tout parahorique $P$, on a une factorisation de $f^{rs}$ en:
\[\wkC_{rs}\stackrel{\pi_{P}^{rs}}{\rightarrow}\wkC_{rs,P}\ra\kC_{rs}\]
où la première flèche est un $W_{P}$-torseur qui fournit donc une famille d'actions compatibles de $(W_P)_{P\in\Par}$ qui produit l'action de Lusztig de $W_{aff}$ sur $\cS_{rs}$ qui coïncide donc avec l'action de $W_{aff}$ sur $\cS_{rs}$ qui provient de \ref{fond-spr}. De plus,  pour tout $\omega\in\Omega$ l'action de $\omega$ entrelace celle de $W_{P}$ et $W_{\omega P\omega^{-1}}$ d'après \eqref{omega}, donc les deux actions de $\widetilde{W}$ sont les mêmes.
\epf

\subsection{$\cP$-équivariance}
Dans le paragraphe précédent, on a vu que le faisceau de Grothendieck-Springer affine était canoniquement muni d'une action de $\widetilde{W}$. 
Nous allons maintenant voir qu'il est également muni d'une action de $\cP$ et que ces actions commutent.
\bprop\label{P-equi}
Le faisceau de Grothendieck-Springer affine $\cS$ est canoniquement $\cP$-équivariant.
\eprop
\bpf
D'après \ref{P-act} et \ref{P-act2}, la flèche $\wkC\ra\kC$ est $\cL G\times\cP$-équivariante et induit donc un morphisme $\cP$-équivariant:
\[[\wkC/\cL G]\ra[\kC/\cL G]\]
Notons pour simplifier $\cX=[\wkC/\cL G]$ et $\cY=[\kC/\cL G]$, alors $\omega_{\cX}$ est $\cP$-équivariant, i.e. on a un isomorphisme canonique
\[\omega_{\cP\times_{\cY}\cX}\cong\act_{\cX}^{!}\omega_{\cX}\stackrel{\sim}{\rightarrow} p_{\cX}^{!}\omega_{\cX}\cong\omega_{\cP\times_{\cY}\cX}.\]
De plus, on a un carré cartésien:
$$\xymatrix{\cP\times_{\cY}\cX\ar[d]_{\ti{f}}\ar@<-.5ex>[r]\ar@<.5ex>[r]&\cX\ar[d]^{f}\\\cP\ar[r]^{p}&\cY},$$
où les flèches horizontales du haut sont respectivement $p_{\cX}$ et $\act$.
On applique $f_!$ et par changement de base propre \ref{Base}, on déduit:
\[\ti{f}_{!}(\act_{\cX}^{!}\omega_{\cX}\stackrel{\sim}{\rightarrow}p_{\cX}^{!}\omega_{\cX})=p^{!}\cS\stackrel{\sim}{\rightarrow}p^{!}\cS,\]
et $\cS$ est $\cP$-équivariant.
\epf

\bthm\label{P-act3}
Le faisceau de Grothendieck-Springer affine $\cS_{\bullet}$ se relève en un faisceau $\widetilde{W}\times\cP$-équivariant.
\ethm
\bpf
On utilise la description en bar-complexe de $[\kC/\cL G]\cong\colim_{[n]\in\Delta^{op}} \kC\times\cL G^{n}$, donc il suffit de montrer  que les pullbacks $\cS_{\bullet}^{(n)}$ de $\cS$ à $\kC\times\cL G^{n}$ pour tout $n$ ont bien des actions de $\widetilde{W}$ et $\cP$ qui commutent entre elles.
Par construction, les actions de $\widetilde{W}$ et $\cP$ commutent avec l'action de $\cL G$ sur $\cS_{\bullet}^{(n)}$. Or, il résulte de \eqref{act-eq}, que l'action de $\cP$ est donnée par la même formule que celle de $\cL G$, donc elle commute à celle de $\widetilde{W}$ comme souhaité.
\epf

\section{\'{E}noncés de Constructibilité}
Le but de cette section est d'établir l'énoncé \ref{W-constr} qui montre que le faisceau de Grothendieck-Springer affine est constructible comme faisceau de $\bql[\widetilde{W}]$-modules.
\subsection{Complexes faiblement constructibles}

\subsubsection{Définitions}\label{loc-fil}

Soit un schéma affine $S$, on définit la sous-catégorie des faiblement constructibles  $\cD_{fc}(X)\subset\cD(X)$, comme la plus petite sous-$\infty$-catégorie pleine, stable par colimites finies et rétractions qui contient les complexes $K\otimes_{\bql}V$ pour $K\in\cD_c(X)$ et $V$ un $\bql$-espace vectoriel.
Pour tout morphisme $f:X\ra Y$, le foncteur $f^!:\cD(Y)\ra\cD(X)$ préserve $\cD_{fc}$. Pour tout préchamp $\cX\in\PrStk_k$, on définit alors:
\[\cD_{fc}(\cX)=\varprojlim\limits_{S\ra\cX}\cD_{f-c}(S),\]
pour $S\in\Aff_k$.
\medskip

Etant donné un $k$-schéma de type fini $S$ et un point géométrique $\bar{s}$, on considère $\Loc(S)$ la catégorie des $\bql$-systèmes locaux constructibles et $\IndLoc(S)$ la ind-catégorie correspondante.
Tout ind-système local $L$ admet une filtration canonique $F_{\bullet}$ telle que $\gr_{F}(L)$ est semisimple: pour tout $i\in\NN$, $\gr_{i}(L)$ consiste en la somme directe des sous-modules simples de $L/L_{\leq i-1}$.

\bdefi
Soit un $k$-schéma de type fini $S$ et $L$ un ind-système local, il est dit faiblement constructible si la filtration $F_{\bullet}$ est finie et pour tout $i$, $gr_{i}(F)$ a un nombre fini de facteurs simples. On note $\FLoc(S)$ la catégorie des ind-systèmes locaux faiblement constructibles.
\edefi
Si $L\in\FLoc(S)$, pour tout $i\in\NN$, $gr_{i}(L)$ est une somme finie  de modules $L_{M}=M\otimes_{\bql}\Hom_{\bql}(M,L)$ pour un sous-système local simple $M\subset L$ de dimension finie. En particulier, comme $gr_{F}(L)$ est finie et comme $\cD_{fc}(S)$ est stable par rétractions et colimites finies on a immédiatement que $\FLoc(S)\subset\cD_{fc}(S)$.

\bprop\label{T-f-c}
Soit $S$ un $k$-schéma de type fini, soit $(S_{\al})$ une stratification finie constructible. On a les assertions suivantes:
\benumr
\remi
 $K\in\cD_{fc}(S)$ si et seulement si pour tout $\al$, $i_{\al}^{!}K\in\cD_{fc}(S_{\al})$ avec $i_{\al}: S_{\al}\ra S$.
\remi
$K\in\cD(S)$ est faiblement constructible si et seulement si il existe une stratification finie constructible $(S_{\al})$, telle que pour tout $i$, $\cH^{i}(i_{\al}^*K)$ et $\cH^{i}(i_{\al}^!K)$ sont des systèmes locaux faiblement constructibles.
\eenum 
\eprop
\bpf
(i) L'implication directe a déjà été vue.  Pour la réciproque, on procède par récurrence sur le nombre de strates et il suffit de considérer le cas où $S=U\coprod F$ avec une immersion fermée $i:F\hra S$ de complémentaire $j:U\ra S$.
Comme $\cD_{fc}$ est stable par colimites finies, en utilisant les suites fibrantes pour $(i^!,j_*)$, il suffit de voir que $i_*$ et $j_*$ préservent $\cD_{fc}$, ce qui est clair puisqu'ils préservent la constructibilité, les rétractions et commutent aux colimites finies.

(ii) Pour l'implication directe, les objets de la forme $K\otimes_{\bql} V$ avec $V\in\Vect_{\bql}$ et $K\in\cD_{c}(S)$ vérifient clairement cette propriété. Il s'agit de voir que cette propriété est stable par colimites finies et rétractions et cela se déduit du fait que la catégorie des ind-systèmes locaux est stable par extensions, noyaux et conoyaux, ce qui est clair. Pour la réciproque cela se déduit de (i) et de l'inclusion $\FLoc(S)\subset\cD_{fc}(S)$.
\epf

\subsection{Rappels sur Barr-Beck-Lurie}
\subsubsection{Modules à gauche sous une algèbre associative}
Soit $\cC$ une $\infty$-catégorie monoïdale symétrique (\cite[Def. 2.0.07]{Lu2}), d'après \cite[Déf.4.1.1.6]{Lu2} on peut considérer $\AlgAss(\cC)$ l'$\infty$-catégorie des algèbres associatives dans $\cC$. Pour $\cA\in\AlgAss(\cC)$, soit $\Mod(\cA)$ l'$\infty$-catégorie des $\cA$-modules en catégories, qui consistent en les $\infty$-catégories qui admettent une action de $\cA$ (\cite[Ch. 0, sect. 3.4]{GR}).
Soit une paire $(\cA,\cM)$ constituée d'une algèbre associative de $\cC$ et d'une $\infty$-catégorie $\cM\in\Mod(\cA)$, on peut considérer la catégorie $\Mod_{\cA}(\cM)$ des objets de $\cM$ qui admettent une structure de $\cA$-module à gauche  (\cite[4.2.1.1.3]{Lu2} ou \cite[Ch. 0, sect. 3.5]{GR}).
Le théorème de Barr-Beck-Lurie consiste à déterminer sous quelles conditions une $\infty$-catégorie se décrit comme une certaine catégorie $\Mod_{\cA}(\cM)$.
\medskip

Si $\cC$ est une $\infty$-catégorie, la catégorie $\Fonct(\cC,\cC)$ admet une structure monoïdale symétrique naturelle de telle sorte que $\cC$ est un module en catégories sous  $\Fonct(\cC,\cC)$ (\cite[Ch. 0, sect. 3.4.6]{GR}). On appelle \textsl{monade}, une algèbre associative $\cA$ de $\Fonct(\cC,\cC)$ agissant sur $\cC$.
Soit alors  $\Mod_{\cA}(\cC)$, la catégorie des objets de $\cC$ avec une structure de $\cA$-module à gauche.
On dispose d'un foncteur d'oubli:
\[\obl_{\cA}:\Mod_{\cA}(\cC)\ra\cC\]
qui admet un adjoint à gauche
\[\ind_{\cA}:\cC\ra\Mod_{\cA}(\cC)\]
qui envoie tout objet $C\in\cC$ sur le $\cA$-module libre à gauche engendré par $C$, de telle sorte que la composée:
\[\obl_{\cA}\circ\ind_{\cA}:\cC\ra\cC\]
s'identifie au foncteur de tensorisation par $\cA$  \cite[4.2.4.8]{Lu2}.
\medskip

\subsubsection{Situation de Barr-Beck-Lurie}
Soient $\cC$ et $\cD$ deux $\infty$-catégories et $\cA\in\AlgAss(\Fonct(\cC,\cC))$, alors d'après \cite[Ch.0, sect. 3.4.6]{GR}, $\Fonct(\cD,\cC)$ admet une structure de module en catégories sur $\Fonct(\cC,\cC)$ et  tout foncteur $G\in\Fonct(\cD,\cC)$ qui admet une structure de $\cA$-module  est équivalent à une factorisation:
\[G:\cD\ra\Mod_{\cA}(\cC)\stackrel{\obl_{\cA}}{\rightarrow}\cC.\]
On  considère alors $G:\cD\ra\cC$ qui admet un adjoint à gauche $F:\cC\ra\cD$. Il résulte de \cite[Ch.0, sect. 3.6.6]{GR} que 
$\cA=G\circ F\in\Fonct(\cC,\cC)$ admet une structure d'algèbre associative et on a une factorisation:
\[G:\cD\stackrel{\tilde{G}}{\rightarrow}\Mod_{\cA}(\cC)\stackrel{\obl_{\cA}}{\rightarrow}\cC.\]
On dit que $G$ est \textsl{monadique} si $\tilde{G}$ est une équivalence.
Le théorème suivant donne un critère de monadicité \cite[4.7.3.5]{Lu2}:
\bthm\label{BBL}[Barr-Beck-Lurie]
Soit $G:\cD\ra\cC$ un foncteur entre $\infty$-catégories avec un adjoint à gauche $F:\cC\ra\cD$. On suppose que:
\benumr
\remi
$G$ est conservatif.
\remi
$G$ préserve les réalisations géométriques, i.e. les colimites indexées sur les ensembles simpliciaux.
\eenum
Alors $G$ est monadique. En particulier, on a une équivalence canonique $\cD\cong\Mod_{\cA}(\cC)$ avec $\cA=G\circ F$.
\ethm

\subsection{Faisceaux $\Gm$-constructibles}
Soit $\Gm$ un groupe discret et $\cX$ un préchamp, $\cD(\cX)$ est une $\infty$-catégorie stable $\bql$-linéaire, donc c'est un module en catégories sur $\Vect_{\bql}$ et $\bql[\Gm]$ est une algèbre associative de $\Vect_{\bql}$.
On définit alors la catégorie des objets $\Gm$-équivariants:
\[\cD^{\Gm}(\cX)=\Mod_{\bql[\Gm]}(\cD(\cX)).\]

Par la suite, on aura besoin de plusieurs descriptions équivalentes de $\cD^{\Gm}(\cX)$.
On considère l'action triviale de $\Gm$ sur $\cX$, soit $p:\cX\ra[\cX/\Gm]$, cette flèche est ind-fp-propre et étale comme $\Gm\times\cX\ra\cX$ l'est.
On dispose d'un foncteur
\[p^{!}:\cD([\cX/\Gm])\ra\cD(\cX).\]
\blem\label{G1}
\benumr
\remi
Le foncteur $p^{!}$ induit une équivalence de catégories $\cD([\cX/\Gm])\stackrel{\sim}{\rightarrow}\cD^{\Gm}(\cX)$.
\remi
De plus, on a également une équivalence canonique :
\begin{equation}
\Mod_{\bql[\Gm]}(\cD(\cX))\cong\Fonct(B\Gm,\cD(\cX)).
\label{GmLu}
\end{equation}
\eenum
\elem

\bpf
(i) Il suffit d'appliquer \ref{BBL}. Comme $p$ est ind-fp-propre, $p^{!}$ admet un adjoint à gauche $p_!$ et  par construction $p_!$ commute aux petites colimites. De plus, $p$ est étale, donc par descente, $p^!$ est conservatif, on peut donc appliquer Barr-Beck-Lurie.
Il suffit donc d'identifier  l'algèbre $\cA$ avec les notations de \ref{BBL}, or pour tout $K\in\cD(\cX)$, on a  d'après \cite[Preuve de 5.6.5]{BKV}:
\[p^{!}p_{!}K\cong K\otimes_{\bql}\bql[\Gm],\]
et le résultat suit.

(ii) En utilisant la description de $[X/\Gm]$ à l'aide du bar-complexe, on a une équivalence de catégories:
\[\cD([\cX/\Gm)\cong\lim_{[n]\in\Delta^{op}}\cD(\Gm^{n}\times_k\cX)\]
Or, comme $\Gm$ est discret, on a $\cD(\Gm^{n}\times_k\cX)\cong\Fonct(\Gm^{n},\cX)$, soit:
\[\lim_{[n]\in\Delta^{op}}\cD(\Gm^{n}\times_k\cX)\cong\Fonct(\colim_{[n]\in\Delta^{op}}\Gm^{n},\cX)\cong\Fonct(B\Gm,\cX),\]
comme souhaité.
\epf

On va utiliser maintenant une troisième description de $\cD^{\Gm}(\cX)$. Pour tout $k$-schéma de type fini $S$ et tout groupe discret $\Gm$, on considère l'$\infty$-catégorie $\cD_{ctf}^{b}(S,\bql[\Gm])$ des complexes bornés constructibles de Tor-dimension finie de $\bql[\Gm]$-modules.
On définit alors:
\[\cD(S,\bql[\Gm])=\Ind(\cD_{ctf}^{b}(S,\bql[\Gm]))\]
 et on étend ensuite par extension de Kan à droite et à gauche à la catégorie de tous les préchamps, les foncteurs $\cD^{b}_{ctf}(-,\bql[\Gm])$ et $\cD(-,\bql[\Gm])$.
Comme un complexe borné constructible de $\bql[\Gm]$-modules, est clairement ind-constructible, à nouveau par fonctorialité de l'extension de Kan, on dispose pour tout préchamp $\cX$ d'un foncteur d'oubli:
\[\cD(\cX,\bql[\Gm])\ra\cD(\cX)\]
qui admet un adjoint à gauche donné par $K\mapsto K\otimes_{\bql}\bql[\Gm]$.
En appliquant \ref{BBL}, on obtient  une équivalence de catégories:
\[\eta:\cD(\cX,\bql[\Gm])\cong\cD^{\Gm}(\cX).\]
En particulier, on peut donc voir la catégorie $\cD^{\Gm}(\cX)$ à la fois, comme la catégorie des faisceaux sur $[\cX/\Gm]$, comme la catégorie des faisceaux de $\bql[\Gm]$-modules sur $\cX$ et comme la catégorie des $\bql$-faisceaux sur $\cX$ qui admettent une structure de $\bql[\Gm]$-module.
\bdefi
Soit $\cX$ un préchamp et un groupe $\Gm$, on définit alors la catégorie des objets $\Gm$-constructibles:
\[\cD^{\Gm}_c(\cX)=\eta(\cD^{b}_{ctf}(\cX,\bql[\Gm])).\]
\edefi
\brem\label{rem-parf}
Il résulte de \cite[Tag. 03TT]{Sta} que pour un $k$-schéma de type fini $X$ et si $\bql[\Gm]$ est noethérien, les complexes de $K\in\cD^{\Gm}_c(\cX)$, consistent en ceux qui peuvent être représentés par un complexe fini de faisceaux de $\bql[\Gm]$-modules, plats et constructibles. En particulier, si $\bql[\Gm]$ est noethérien, les fibres $K_{\bar{x}}$ en tout point géométrique $\bar{x}\in X$  sont des complexes parfaits de $\bql[\Gm]$-modules. Si de plus $\bql[\Gm]$ est de dimension cohomologique finie, alors les complexes parfaits de $\bql[\Gm]$-modules, sont simplement les complexes de $\bql[\Gm]$-modules dont les faisceaux de cohomologie sont des $\bql[\Gm]$-modules de type fini.
\erem
\blem\label{Get-loc}
Les foncteurs $\cD^{\Gm}_c$ et $\cD^{\Gm}$ sont des faisceaux pour la topologie étale.
\elem
\bpf
En vertu du lemme \ref{G1} et comme $\cD$ est un faisceau pour la topologie étale, on a immédiatement que $\cD^{\Gm}$ est aussi un faisceau pour la topologie étale. Ainsi, pour établir que $\cD^{\Gm}_{c}$ est un faisceau étale, il s'agit de vérifier que pour tout morphisme étale surjectif $h:\cY\ra\cX$ de préchamps, si $h^{!}K\in\cD^{\Gm}_{c}(\cY)$, alors $K\in\cD^{\Gm}_{c}(\cX)$. Par définition du foncteur $\cD^{\Gm}_{c}$, on se ramène immédiatement au cas où $\cX=X$ est un schéma affine et quitte à raffiner le recouvrement, on peut supposer $Y$ affine. On choisit alors une présentation $X\simeq\varprojlim X_{\al}$ où les $X_{\al}$ sont des $k$-schémas affines de type fini, ainsi $h$ se descend en un morphisme étale surjectif de schémas affines $h_{\al}: Y_{\al}\ra X_{\al}$  et $Y\simeq\varprojlim_{\beta\geq\al}Y'_{\beta}=X_{\beta}\times_{X_{\al}} Y_{\al}$.
Comme $\cD^{\Gm}_{c}(X)\simeq\colim\cD^{\Gm}_{c}(X_{\al})$ et pareillement pour $\cD^{\Gm}_c(Y)$, on est ramené au cas où $h$ est étale entre $k$-schémas  affines de type fini et là c'est clair.
\epf
\subsubsection{Fonctorialité}
Pour tout morphisme de $k$-schémas de type fini $f:X\ra Y$, $f^!$ préserve $\cD^{b}_{ctf}()$ d'après \cite[Exp. 7, 1.7]{SGA4.5}, il en est donc de même, par extension de Kan, pour tout morphisme de préchamps $f:\cX\ra\cY$ et donc $f^{!}$ préserve en particulier les objets $\Gm$-constructibles.
De plus, si $X$ est de type fini et $X_{\al}\stackrel{i_{\al}}{\rightarrow}X$ une stratification finie constructible, on obtient alors que:
\[K\in\cD^{\Gm}_{c}(X)\Longleftrightarrow\forall~\al, i_{\al}^{!}K\in\cD^{\Gm}_{c}(X_{\al}).\]
L'implication directe vient de ci-dessus et la réciproque vient de la stabilité par les opérations de $\cD^{b}_{ctf}$ en considérant les suites exactes fibrantes successives et stabilité par passage au cône de $\cD^{b}_{ctf}$.
Comme un faisceau constructible de $\bql[\Gm]$-modules est clairement faiblement constructible, on obtient que le foncteur d'oubli:
\begin{equation}
\obl:\cD^{\Gm}(\cX)\ra\cD(\cX)
\label{obl-fc}
\end{equation}
envoie $\cD^{\Gm}_c(\cX)$  dans $\cD_{fc}(\cX)$ pour tout préchamp $\cX$.
\bdefi
Pour tout groupe discret $\Gm$ et un préchamp $\cX$, on note $\cD_{fc}^{\Gm}(\cX)$ l'ensemble des $K\in\cD^{\Gm}(\cX)$ tels que $\obl(K)\in\cD_{fc}(\cX)$.
\edefi
\medskip

Pour tout morphisme $\al:\Gm_1\ra\Gm_2$, on dispose d'un foncteur d'induction:
\begin{equation}
\ind^{\Gm_2}_{\Gm_1}:\cD(\cX,\bql[\Gm_1])\ra\cD(\cX,\bql[\Gm_2])
\label{def-ind}
\end{equation}
donné par $K\mapsto K\otimes^{L}_{\bql[\Gm_1]}\bql[\Gm_2]$,
qui préserve $\cD_{ctf}^{b}$ et commute à $f^!$ pour tout morphisme de préchamps $f:\cY\ra\cX$. Si l'on considère le cas particulier de la projection $\Gm\ra\{1\}$ alors $\ind^{\Gm}_{\{1\}}$ est le foncteur des coinvariants $\coinv_{\Gm}$.
Enfin, pour toute représentation de dimension finie $\tau\in\Rep^{df}_{\bql}(\Gm)$ et $K\in\cD^{\Gm}_{c}(\cX)$, $K\otimes_{\bql}\tau\in\cD^{\Gm}_{c}(\cX)$. On peut donc définir la composante $\tau$-isotypique de $K$ par:
\[K_{\tau}=\coinv_{\Gm}(K\otimes_{\bql}\tau)\]
Il résulte de ci-dessus que:
\begin{equation}
\forall~ \tau\in\Rep^{df}_{\bql}(\Gm),\forall~ K\in\cD^{\Gm}_{c}(\cX), K_{\tau}\in\cD^{\Gm}_c(\cX).
\label{tcoinv}
\end{equation}
\blem\label{res-cons}
Soit $\cX$ un $\infty$-champ qui satisfait le recollement, $(\cX_{\al})$ une stratification finie constructible, alors on a:
 \[
K\in\cD^{\Gm}_{c}(\cX)\Longleftrightarrow\forall~\al, i_{\al}^!\cX_{\al}\in\cD^{\Gm}_{c}(\cX_{\al}),\]
pour $i_\al:\cX_{\al}\ra\cX$.
\elem
\bpf
Le sens direct est clair par stabilité de $\cD^{\Gm}_{c}$ par $f^!$. Pour la réciproque, si $K\in\cD^{\Gm}(\cX)$ tel que  pour tout $\al$, $i_{\al}^!\cX_{\al}\in\cD^{\Gm}_{c}(\cX_{\al})$, en considèrant les suites fibrantes et en raisonnant par récurrence, on se ramène à la situation d'un fermé de présentation  finie $i:\cX_{F}\hra\cX$ de complémentaire $j:\cX_{U}\hra\cX$ et il s'agit de montrer que :
\[i_*i^{!}K\in\cD^{\Gm}_{c}(\cX), j_*j^{*}K\in\cD^{\Gm}_{c}(\cX),\]
et cela se ramène alors au fait que $j_*$ et $i_*$ préservent les constructibles comme $\cX$ satisfait le recollement d'après \ref{glustk}(iv)\footnote{A strictement parler, il est seulement établi que $j_*$ et $i_*$ préservent les $\bql$-constructibles et pas les $\bql[\Gm]$-constructibles, mais l'argument s'étend une fois que l'on a la stabilité des $\bql[\Gm]$-constructibles pour les $k$-schémas de type fini, ce qui est le cas.}.
\epf
On a l'énoncé de descente suivant:
\bprop\label{lis-desc}
\benumr
\remi
Soit $f:\cX\ra\cY$ un morphisme lisse surjectif de préchamps au sens de \ref{fpl}, alors pour tout $K\in\cD^{\Gm}(\cY)$, on a $f^{!}K\in\cD^{\Gm}_c(\cX)$ si et seulement si $K\in\cD^{\Gm}_c(\cY)$.
\remi
Soit $f:X\ra Y$ un morphisme de présentation finie, surjectif de schémas placides, on suppose $\bql[\Gm]$ noethérien, de dimension cohomologique finie. Soit $K\in\cD^{\Gm}(Y)$, alors on a $f^{!}K\in\cD^{\Gm}_c(X)$ si et seulement si $K\in\cD^{\Gm}_c(Y)$.
\eenum
\eprop
\brem
On rappelle que notre notion de lissité n'impose pas d'être localement de présentation finie. L'hypothèse que $\bql[\Gm]$ est de dimension cohomologique finie n'est probablement pas nécessaire, mais suffisant pour nos applications et raccourcit les arguments.
\erem
\bpf
(i) Il suffit de vérifier l'énoncé après tiré-en-arrière sur tout schéma affine. On peut donc supposer que $\cY=Y$ est affine et $\cX=X$ est un schéma. De plus, quitte à recouvrir $X$ par des ouverts affines et comme $Y$ est quasi-compact, on peut supposer $X$ affine également.
Il résulte alors de la description \ref{G1} et du fait que $\cD^{\Gm}_c$ est local pour la topologie étale, que l'on peut se ramener au cas où 
$X\ra Y$ est fortement pro-lisse et en considérant une présentation $X\simeq\varprojlim X_{\al}$, on a alors:
\[\cD_c^{\Gm}(X)\simeq\colim\cD_c^{\Gm}(X_{\al})\]
et on se ramène au cas où $f$ est lisse. Dans ce cas, $f$ admet localement des sections pour la topologie étale et  le résultat suit en utilisant \ref{Get-loc}.

(ii) A nouveau, comme $\cD^{\Gm}_c$ est un faisceau pour la topologie étale, on peut supposer que $Y$ est affine placidement présenté.
On considère alors une présentation placide $Y\simeq\varprojlim Y_{\al}$, de telle sorte que par définition, on a $\cD^{\Gm}_c(Y)\simeq\colim\cD^{\Gm}_{c}(Y_{\al})$ et le morphisme $f$ se descend en un morphisme surjectif de présentation finie $f_{\al}:X_{\al}\ra Y_{\al}$. On est donc ramené au cas d'un morphisme surjectif entre $k$-schémas de type fini. Par dualité de Verdier, on obtient donc que $f^{*}K\in\cD_{ctf}(X,\bql[\Gm])$. Comme $\bql[\Gm]$ est noethérien de dimension cohomologique finie, on a $\cD_{ctf}(X,\bql[\Gm])=\cD_{c}(X,\bql[\Gm])$. En passant aux $\cH^{i}$, on se ramène alors à un faisceau constructible de $\bql[\Gm]$-modules et cela se déduit de \cite[Tag. 095Q]{Sta}.
\epf
L'énoncé suivant fournit un critère de $\Gm$-constructibilité.
\bprop\label{gm-crit}
Soit $\cX$ un champ placide et $K\in\cD_{fc}^{\Gm}(\cX)$, on suppose $\bql[\Gm]$ noethérien. Alors $K\in\cD^{\Gm}_{c}(\cX)$ si et seulement si pour tout $x\in\cX$, $i_x^{!}K\in\cD^{\Gm}_{c}(\pt)$.
Si de plus, $\bql[\Gm]$ est de dimension cohomologique finie, il suffit de vérifier que pour tout $x$, $i_x^{!}K$ est un complexe borné de $\bql[\Gm]$-modules de type fini.
\eprop
\bpf
Si $K\in\cD^{\Gm}_{c}(\cX)$ alors pour tout morphisme de préchamps, $f^!$ préserve les $\Gm$-constructibles d'après ci-dessus, donc le sens direct est clair.
Réciproquement, en utilisant \ref{lis-desc} on se ramène au cas où $\cX=X$ est affine placidement présenté.
En considérant une présentation placide $X\simeq\varprojlim X_{\al}$ de telle sorte que $\cD^{\Gm}_c(X)\simeq\colim\cD_c^{\Gm}(X_{\al})$ et pareillement pour $K\in\cD_{fc}^{\Gm}(\cX)$, on peut supposer que $X$ est un $k$-schéma de type fini.
Soit $K\in\cD_{fc}^{\Gm}(X)$, il suffit de voir que pour tout $i$, $H^{i}(K)\in\cD^{\Gm}_{c}(K)$. Comme  $\bql[\Gm]$ noethérien, d'après \cite[Tag.03TT]{Sta}, il est d'amplitude finie et on suppose donc que $K$ est un faisceau ind-contructible de $\bql[\Gm]$-modules. D'après \ref{T-f-c}, il existe alors une stratification $(S_\al)$ tel que $!$-restriction de $K$ à $S_{\al}$ est un ind-système local faiblement constructible et d'après \ref{res-cons}, on peut se ramener à $\obl(K)$ qui est un ind-système local faiblement constructible.
Maintenant la filtration \ref{loc-fil} de $\obl(K)$ est $\Gm$-équivariante (on peut également la définir directement pour les systèmes locaux de $\bql[\Gm]$-modules) et est finie comme $\obl(K)$ est faiblement constructible et comme $\bql[\Gm]$ est noethérien, on peut supposer que $\obl(K)$ semisimple, $\Gm$-constructible et  la décomposition en sous-modules simples est aussi $\Gm$-équivariante:
\[\obl(K)=\bigoplus (M\otimes_{\bql} Hom(M,\obl(K)).\]
On est ramené alors au cas où $\obl(K)=M\otimes_{\bql} \Hom(M,\obl(K))$ et comme $K$ est $\Gm$-équivariant, $\Hom(M,\obl(K))$ est une $\Gm$-représentation. De plus pour tout $x\in X$, on a :
\[i_x^{!}K=i_x^{!}M\otimes_{\bql}\Hom(M,\obl(K))\in\cD_{c}^{\Gm}(\pt),\]
 donc $\Hom(M,\obl(K))$ est un complexe parfait de $\bql[\Gm]$-modules et le résultat suit.
Si de plus, $\bql[\Gm]$ est de dimension cohomologique finie, cela se déduit  de \ref{rem-parf}.
\epf

\bthm\label{W-constr}
Le faisceau $\cS_{\bullet}\in\cD^{\widetilde{W}}([\kC_{\bullet}/\cL G])$ est $\widetilde{W}$-constructible. En particulier pour toute représentation de dimension finie $\tau\in\Rep_{\bql}(\widetilde{W})$, la composante $\tau$-isotypique $\cS_{\tau}$ est dans $\cD_{c}([\kC_{\bullet}/\cL G])$.
\ethm
\bpf
D'après \ref{res-cons}, on peut se ramener à montrer la $\widetilde{W}$-constructibilité après restriction aux strates $[\kC_{w,\br}/\cL G]_{red}$. On utilise alors \ref{lis-desc} pour montrer la $\widetilde{W}$-constructibilité après tiré-en-arrière à $\kt_{w,\br}$.

Rappelons le diagramme suivant \eqref{wr-diag}:
\begin{equation}
\xymatrix{\ti{X}_{w,r}\ar[d]_{\ti{f}_{w,r}}\ar[r]&[\wkC_{w,\br}/\cL G]\ar[r]\ar[d]^{\ov{p}}&[\wkC_{\bullet}/\cL G]\ar[d]\\\kt_{w,r}\ar[r]^-{\psi_{w,\br}}&[\kC_{w,\br}/\cL G]\ar[r]&[\kC_{\bullet}/\cL G]}.
\end{equation}
de telle sorte que $f_{w,r}:X_{w,r}=\ti{X}_{w,r,red}\ra\kt_{w,r}$ est localement de présentation finie avec une action d'un réseau $\La_w$ et l'on a une factorisation:
\[X_{w,r}\stackrel{\pi}{\rightarrow}[X_{w,r}/\La_w]\stackrel{\mu}{\rightarrow}\kt_{w,r},\]
où la deuxième flèche est fp-propre surjective et la première étale de groupe $\La_w$.
Posons $\cF_{w,\br}=(f_{w,r})_{!}\omega_{X_{w,r}}=\mu_{!}\pi_{!}\omega_{X_{w,r}}.$
Comme $\pi$ est étale de groupe $\La_w$, on a  d'après  \cite[5.6.5]{BKV}, $\pi^{!}\pi_{!}\omega_{X_{w,r}}\cong\omega_{X_{w,r}}\otimes_{\bql}\bql[\La_{w}]$.
Ainsi, $\pi^{!}\pi_{!}\omega_{X_{w,r}}$ est $\La_w$-constructible donc d'après \ref{lis-desc}, $\pi_{!}\omega_{X_{w,r}}$ est $\La_w$-constructible.
Et d'après \cite[Exp. 7, Thm.1]{SGA4.5}, $\mu$ préserve les $\La_w$-constructibles\footnote{Comme $\mu$ est fp-propre, par descente noethérienne, on commence par se ramener au cas des $k$-schémas de type fini et ensuite on se réduit au cas des coefficients de torsion}.
On obtient donc que $\cF_{w,r}$ est $\La_w$-constructible, donc en particulier essentiellement constructible.
Maintenant $\bql[\widetilde{W}]$ est noethérien de dimension cohomologique finie et comme l'action de $\widetilde{W}$ est la même que celle de Lusztig, d'après \ref{lusz}, d'après Yun \cite[Cor. 4.4.5]{BeKV}, $H^{a}(i_{x}^!\cS)$ est un $\bql[\widetilde{W}]$-module de type fini, donc on conclut par \ref{gm-crit}.
Il ne reste plus qu'à appliquer \eqref{tcoinv} pour en déduire l'énoncé sur les coinvariants.
\epf
\section{Perversité des coinvariants}

\subsection{Un lemme d'homotopie}
La preuve de la proposition suivante nous a été suggérée par Y. Varshavsky.
\bprop\label{homo1}
Soit $P$ un schéma en groupes lisse à fibres géométriquement connexes sur un $k$-schéma de type fini $S$.
Considérons le foncteur de restriction:
\[\res_{P}:\cD_{P}(S)\ra\cD_{P(S)}(S),\]
alors on a un relèvement canonique :
$$\xymatrix{&\cD(S)\ar[d]^{\triv}\\\cD_{P}(S)\ar@{.>}[ur]^\phi\ar[r]^{\res_{P}}&\cD_{P(S)}(S).}$$
\eprop
\bpf
Pour tout schéma en groupes lisse $P$ sur $S$, en utilisant Barr-Beck-Lurie \ref{BBL}, on a une équivalence de catégories:
\[\cD_{P}(S)\cong H_{*}(P/S)-Mod,\]
où l'on voit l'homologie $H_*(P/S)=p_{!}p^{!}\mathbb{Q}_{\ell,S}$ comme un faisceau en DG-algèbres.
En particulier, la proposition se ramène à montrer que la flèche de DG-algèbres:
\[H_*(P(S)/S))\ra H_*(P/S)\]
se factorise par la section unité $1_{S}:\mathbb{Q}_{\ell,S}\ra H_{*}(P/S)$.
Tout, d'abord, notons que $H_{*}(P/S)\in D^{\leq 0}(S)$ pour la $t$-structure standard et que comme $P$ est lisse à fibres géométriquement connexes et que la formation de $H_{*}(P/S)$ commute aux $*$-pullbacks, la composée:
\begin{equation}
\mathbb{Q}_{\ell,S}\stackrel{1_{S}}{\rightarrow} H_{*}(P/S)\ra H_{0}(P/S),
\label{UnS}
\end{equation}
est un isomorphisme.
On déduit ainsi de \eqref{UnS} que la composée:
\[H_{*}(P(S)/S)\ra H_{*}(P)\ra H_{0}(P/S),\]
se relève uniquement en un morphisme $H_{*}(P(S)/S)\ra\mathbb{Q}_{\ell,S}$.
De plus, Comme $P(S)$ est discret, le  faisceau d'algèbres $H_{*}(P(S)/S)$  est concentré en degré zéro et comme $H_{*}(P/S)\in D^{\leq 0}(S)$, on en déduit que la flèche:
\[\Hom(H_{*}(P(S)/S),H_{*}(P/S))\ra\Hom(H_{*}(P/S), H_{0}(P/S))\]
est une équivalence, de telle sorte que le morphisme $H_{*}(P(S)/S)\ra\mathbb{Q}_{\ell,S}$ relève l'application $H_*(P(S)/S))\ra H_*(P/S)$, comme souhaité.
\epf

\bthm\label{cor-homo}
Soit un $k$-schéma de type fini $S$, $P$ un $S$-schéma en groupes lisse. On considère $K\in\cD(S)$ un faisceau $P$-équivariant, alors l'action se factorise par le faisceau des composantes connexes relatives $\pi_{0}(P/S)$, i.e. pour tout $\phi:U\ra S$ l'action de $P(U)$ sur $\phi^!K$ se factorise par $\pi_0(P/S)(U)$.
\ethm
\brem
On renvoie à \cite[Prop. 6.2]{N2}, pour la définition du faisceau $\pi_{0}(P/S)$.
\erem
\bpf
Soit $U$ un $S$-schéma, soit $P^{0}$ la composante neutre de $P$; il existe alors $U'\ra U$ un recouvrement étale tel que l'on a
$\pi_{0}(P)(U')=P(U')/P^{0}(U')$. 
De plus, on a un carré cartésien de champs:
$$\xymatrix{BP^{0}(U')\ar[d]\ar[r]\ar[d]&\pt\ar[d]\\BP(U')\ar[r]&B\pi_{0}(P)(U')}$$
qui induit en vertu de \eqref{GmLu} un diagramme cartésien de catégories:
$$\xymatrix{\cD_{\pi_{0}(P)(U')}(U')\ar[d]\ar[r]&\cD_{P(U')}(U')\ar[d]\\\cD(U')\ar[r]&\cD_{P^{0}(U')}(U')}$$
Le lemme \ref{homo1} appliqué à $P^0$ donne que la composée:
\[\cD_{P}(S)\ra\cD_{P^0}(U')\ra\cD_{P^{0}(U')}(U')\]
se factorise canoniquement par $\cD(U')\ra\cD_{P^{0}(U')}(U')$, on obtient donc une flèche canonique:
\[\cD_{P}(S)\ra\cD_{\pi_{0}(P)(U')}(U'),\]
compatible à la localisation.
Il reste donc à voir que le foncteur $U\mapsto\cD_{\pi_{0}(P)(U)}(U)=\cD(B\pi_{0}(P(U))$ est un faisceau étale en catégories.
Pour simplifier les notations, notons $\Pi=\pi_{0}(P/S)$, on a alors:
\[\cD_{\Pi(U)}(U)=\lim_{[m]\in\Delta^{op}}\cD(\Pi(U)^{m}\times U)=\lim_{[m]\in\Delta^{op}}\cD(\Pi^{m}(U)\times U).\]
Comme la propriété d'être un faisceau étale est stable par limites (\cite[Cor. 5.5.2.4]{Lu1}), il suffit de prouver que le préfaisceau $U\mapsto\cD(\Pi^m(U)\times U)$ est un faisceau étale, ce qui est alors clair puisque $\Pi^m$ et $\cD$ sont des faisceaux étales.
\epf
On en déduit le corollaire suivant pour la fibration de Springer affine.
\bcor\label{Pic-pi}
Pour tout point géométrique $\g\in\kC_{\bullet}(k)$, $a=\chi(\g)$, l'action de $\cP_a$ sur $i_{\g}^!\cS\cong R\Gm(X_{\g},\omega_{X_{\g}})$ se factorise par $\pi_0(\cP_a)$.
\ecor
\brem
On insiste sur le fait que lorsque l'on dit que l'action se factorise, c'est au sens de \ref{cor-homo}.
\erem
Pour tout $\g\in\kC_{\bullet}(k)$ avec $a=\chi(\g)$, Yun définit \cite[sect. 2.7]{Yun} un morphisme:
\begin{equation}
\sigma_\g:\bql[X_*(T)]^W\ra \bql[\pi_0(\cP_a)].
\label{sigma}
\end{equation}
On a la conjecture locale qui est une généralisation du théorème de Yun \cite[Thms. 1-2]{Yun}.
\bconj\label{act-sph}
Pour tout corps algébriquement clos $K$ et $\g\in\kC_{\bullet}(K)$, l'action de la partie sphérique $\bql[X_*(T)]^W$ de $\widetilde{W}$ sur $R\Gm(X_{\g},\omega_{X_{\g}})$ se factorise par \eqref{sigma}.
\econj
\brem
On peut également formuler une variante pour la cohomologie à support compact, mais nous n'en aurons pas besoin dans ce travail. L'avantage de travailler avec des complexes plutôt qu'avec des faisceaux de cohomologie est que d'après \cite[5.8]{YunII}, on a une dualité entre l'homologie et la cohomologie à support compact, de telle sorte que les deux énoncés devraient être équivalents.
\erem
Dans la suite, il sera commode de considérer un relèvement de $\sigma_{\g}$ en un morphisme $X_{*}(T)\ra\pi_0(P_a)$.
D'après Yun \cite[2.7]{Yun}, considérons le revêtement caméral $B_{a}=k[[t]]\otimes_{\cO_{\kc}}\cO_{\kt}$ obtenu via la flèche $a:\Spec(k[[t]])\ra\kc$. Si l'on choisit une composante $B$ de la normalisation $B_{a}^{\flat}$ de $B_a$, alors on obtient une surjection:
\begin{equation}
\sigma'_{\g}:X_*(T)\ra\pi_{0}(P_a)
\label{w-lift}
\end{equation}
et si l'on modifie le choix de la composante, le morphisme $\sigma'_{\g}$ diffère par l'action de $W$ sur $X_{*}(T)$.

\subsection{Application à la fibration de Hitchin}
Le lemme d'homotopie permet de formuler des conjectures locales et globales sur l'action du Picard sur la fibration de Grothendieck-Springer affine et sur la fibration de Hitchin.

Pour toute courbe projective lisse $X$ sur un corps $k$ et un diviseur effectif $D$ sur $X$, on dispose de la fibration de Hitchin parabolique $\cM^{par}$ qui classifie les quadruplets $(x,E,\phi, E_{B})$ où $x\in X$, $E$ est un $G$-torseur sur $X$, $\phi\in H^{0}(X,\ad(E)\otimes_{k} \cO_{X}(D))$ où $\ad(E)$ est le fibré adjoint de $E$ et $E_{B}$ est un $B$-réduction de $E$ en $x$, compatible à $\phi$.
On a un morphisme de Hitchin parabolique, induite par le polynôme caractéristique:
\[f^{par}:\cM^{par}\ra X\times\cA\]
avec $\cA=H^{0}(X,\kc\otimes_{k}\cO(D))$, donnée par $(x,E,\phi, E_{B})\mapsto (x,\chi(\phi))$.
Yun construit dans (\cite{Yun},\cite{YunII}) une action de $\widetilde{W}$ sur le complexe $f^{par}_*\bql$ au-dessus de l'ouvert génériquement régulier semisimple $\cA^{\heartsuit}\subset\cA$. La fibration de Hitchin parabolique admet également l'action d'un champ de Picard global $\cP$ au-dessus de $X\times\cA$ qui induit une action, commutant à celle de $\widetilde{W}$, sur $f^{par}_*\bql$.
On a également un morphisme global de faisceaux au-dessus de $\cA$:
\begin{equation}
\sigma:\bql[X_*(T)]^{W}\ra\bql[\pi_{0}(\cP/\cA)],
\label{sph2}
\end{equation}
où $\pi_{0}(\cP/\cA)$ est le faisceau des composantes connexes relatives.
Dans \cite[Thm.3]{YunII}, Yun montre que si l'on restreint l'action de $\widetilde{W}\times\cP$ à $\bql[X_*(T)]^{W}\times\cP$ sur $f^{par}_*\bql$ et qu'on la semisimplifie, i.e. on regarde l'action induite sur les faisceaux de cohomologie pervers $R^{i}f^{par}_*\bql$ alors elle se factorise en un action diagonale de $\pi_{0}(\cP/\cA)\times \pi_{0}(\cP/\cA)$ via \eqref{sph2}.
Une des raisons pour laquelle Yun se restreint à l'action semisimplifiée est qu'il ne disposait du lemme d'homotopie que sur les $R^{i}f^{par}_*\bql$. En appliquant \ref{homo1}, on obtient immédiatement donc l'énoncé suivant:
\bthm\label{fact}
L'action de $\cP$ sur le complexe $f^{par}_*\bql$ se factorise en une action de $\pi_{0}(\cP/\cA)$.
\ethm
Cela permet donc de formuler la variante globale de la conjecture \ref{act-sph}.

\bconj\label{act-sph2}
La restriction de l'action de $\widetilde{W}$ à $\bql[X_*(T)]^{W}$  sur $f^{par}_*\bql$ se factorise par $\pi_{0}(\cP/\cA)$.
\econj
\brems
\remi
En fait, une inspection détaillée de la preuve de Yun dans le cas global montre qu'une fois que l'on a le lemme d'homotopie \ref{homo1}, les arguments de \cite{Yun} s'étendent mutatis mutandis pour déduire \ref{act-sph2}. En revanche, pour la conjecture locale, il semble plus délicat de le déduire de la conjecture globale.
\remi
On remarquera la dissymétrie entre la conjecture locale qui se formule point par point alors que la conjecture globale se formule en famille, de manière faisceautique. Cela vient du fait que localement, il est délicat de définir le faisceau des composantes connexes relatives $\pi_{0}(\cP/\clp\kc)$ ainsi que l'analogue local du morphisme $\sigma$ qui devrait être défini au-dessus de $\clp\kc$.
La construction du morphisme global nécessite d'utiliser le revêtement caméral universel $\ti{\cA}\ra X\times\cA$; on construit d'abord la flèche $\sigma$ au-dessus de $X\times \cA$, puis on montre que la flèche se descend au-dessus de $\cA$. L'analogue local de ce revêtement caméral vit sur $\Spec(k[[t]]])\hat{\times}\clp\kc$ et on se retrouve confronté à la difficulté d'avoir à considérer des anneaux de séries formels $A[[t]]$ pour des anneaux non-noethériens et ensuite de montrer que la flèche descend le long du morphisme fpqc $\Spec(A[[t]])\ra\Spec(A)$.
\erems
\subsection{Calcul des coinvariants}
Soit $\tau$ une représentation de dimension finie de $\widetilde{W}$, on dit qu'elle est de torsion si sa restriction $\tau_{\vert X_{*}(T)}$ se décompose en une somme directe de caractères de torsion.
L'énoncé principal de la section est le suivant:
\bthm\label{coinv}
Soit $G$ semisimple simplement connexe, soit $\tau\in\Rep_{\bql}^{df}$ de torsion et le faisceau $S_{\tau}=\coinv_{\tilde{W}}(\cS_{\bullet}\otimes_{\bql}\tau)$ des $\tau$-coinvariants. On suppose \ref{act-sph} vérifiée, alors $\cS_{\tau}$ est pervers.
\ethm
\brem
A priori, ces faisceaux pervers ont en revanche des supports. La preuve du théorème nous donne une borne supérieure que l'on explicite dans \ref{G-supp}.
\erem
On commence par quelques dévissages.

\subsubsection{Coinvariants sous $\widetilde{W}$}
Soit $G$ semisimple simplement connexe.
Tout d'abord, il résulte de \cite[Rmq., Thm. 2.16]{Mi} que l'on a une équivalence:
\[\pt\cong\colim_{I\subset P} \widetilde{W}/W_{P},\]
où $P$ parcourt les parahoriques standards de $\cL G$. En quotientant par $\widetilde{W}$, on en déduit une équivalence au niveau des classifiants:
\begin{equation}
B\widetilde{W}\cong\colim_{I\subset P} BW_{P}.
\label{aff-eq}
\end{equation}
Maintenant d'après \eqref{GmLu} pour tout préchamp, on a une équivalence $\cD^{\widetilde{W}}(\cX)\cong\Fonct(B\widetilde{W},\cD(\cX))$, de telle sorte que l'on en déduit une équivalence:
\[\cD^{\widetilde{W}}(\cX)\cong\Fonct(\colim BW_P,\cD(\cX))\cong\varprojlim_{I\subset P}\cD^{W_{P}}(\cX)\]
et l'on obtient ainsi une équivalence :
\begin{equation}
\cS_{\tau}\cong\colim_{I\subset P}\cS_{\tau,P}.
\label{W-coinv}
\end{equation}
où $\cS_{\tau,P}$ les $\tau$-coinvariants pris relativement à $W_P$.

\blem\label{devis-1}
Pour établir \ref{coinv}, il suffit de  montrer que pour toute strate $(w,r)$, on a:
\begin{equation}
 i_{w,r}^{!}\cS_{\tau}\in\mathstrut^{p}\cD^{\geq-\nu_{w,r}}(\kC_{w,r}).
\label{p-cond3}
\end{equation}
De plus, l'inégalité est stricte pour toute strate $\kC_{w,r}$ pour laquelle $d(w,r)\geq 1$.
\elem
\bpf
D'après \ref{fond-spr}, $\cS$ est pervers. Pour un groupe fini , les $\tau$-coinvariants sont $t$-exacts  et d'après \cite[Lem. 6.1.2]{BKV}, $\cD^{\leq 0}$ est stable par toutes les petites colimites, de telle sorte qu'il résulte de \eqref{W-coinv} que le foncteur des $\tau$-coinvariants est $t$-exact à droite. Il suffit donc de vérifier la moitié des inégalités \eqref{t1} et \eqref{t2} et le lemme suit.
\epf
On va se ramener maintenant à un énoncé fibre à fibre.
Tout d'abord, on rappelle que l'on a une flèche lisse:
\[\psi_{w,r}:\kt_{w,r}\ra[\kC_{w,r}/\cL G]_{red}\]
de telle sorte que d'après \cite[Prop. 6.3.3]{BKV}, $\psi_{w,r}^!$ est $t$-exact. Il suffit donc de vérifier l'énoncé après tiré-en-arrière à $\kt_{w,r}$.

On a alors le lemme suivant :
\blem\label{t-pt}
Soit un schéma affine $X$, placidement présenté et équidimensionnel, muni de la $t$-structure de \ref{t-codim}, alors on a $K\in\mathstrut^{p}\cD^{\geq 0}(X)$ si et seulement si pour tout $x\in X$, on a $i_x^!K\in\mathstrut^{p}\cD^{\geq 0}(\pt)$.
\elem
\bpf
On commence par choisir une présentation placide $X\simeq\varprojlim X_{\al}$. Il résulte alors de \cite[Prop. 6.1.4]{BKV} que l'on a :
\[\mathstrut^{p}\cD^{\geq 0}(X)\simeq\colim\mathstrut^{p}\cD^{\geq 0}(X_{\al})\]
et que les flèches de transition sont $t$-exactes.
On se ramène ainsi immédiatement au cas où $X$ est un $k$-schéma de type fini équidimensionnel.
Comme $\cD(X)=\Ind(\cD_c)$, il résulte de \cite[Lem. 6.1.2]{BKV}, que l'on peut supposer $K\in\cD_c(X)$. Dans ce cas, la $t$-structure est juste le décalage par la dimension de $X$ de la $t$-structure perverse standard et c'est bien connu \cite[2.2.12.(ii)]{BBD}.
\epf
Ainsi, en vertu de \ref{t-pt} et comme $\psi_{w,r}$ est $t$-exact, il suffit donc de montrer l'inégalité suivante:
\begin{equation}
\forall x\in\kt_{w,r}, i_x^{!}\cS_{\tau}\in\mathstrut^{p}\cD^{\geq-\nu_{w,r}}(\pt),
\label{p-fib}
\end{equation}
où $i_x$ désigne la composée $x\ra\kt_{w,r}\ra[\kC_{\bullet}/\cL G]$.
\subsubsection{Dévissage des $\widetilde{W}$-coinvariants}
Dans la suite, tous les produits tensoriels sont dérivés.
Soit un corps algébriquement clos $K$ et $i_{\g}:\Spec(K)\ra\kC_{w,r}$, le changement de base propre  donne alors:
\[i_\g^{!}\cS_{\tau}\cong (R\Gm(X_{\g},\omega_{X_{\g}})\otimes_{\bql}\tau)\otimes_{\bql[\widetilde{W}]}\bql.\]
 On a une décomposition en produit semi-direct $\widetilde{W}=X_*(T)\rtimes W$.
On obtient ainsi l'identité:
\[\coinv_{\widetilde{W}}=\coinv_{W}\circ\coinv_{X_*(T)}.\]
Il nous suffit donc de vérifier \eqref{p-fib} pour les coinvariants sous $X_*(T)$.
Par hypothèse, $\tau_{\vert X_*(T)}$ est une somme directe de caractères de torsion, on se ramène donc au cas où $\tau_{\vert X_*(T)}=\kappa:X_*(T)\ra\bZ/N\bZ$ pour un certain $N\in\NN$. On note $\bqlk$ la représentation de dimension un associée.
On va voir qu'il suffit de calculer ces coinvariants lorsque l'on se restreint à $\bql[X_*(T)]^{W}$.

Comme $G$ est simplement connexe, d'après \cite[Thm.6.1.2]{CG},  $\bql[X_*(T)]$ est un $\bql[X_*(T)]^{W}$-module libre de type fini et on a un isomorphisme canonique $W$-équivariant de $\bql[X_*(T)]^{W}$-modules :
\[\bql[X_*(T)]\cong\bql[X_*(T)]^{W}\otimes_{\bql}\bql[W].\]
On obtient ainsi que $R\Gm(X_{\g},\omega_{X_{\g}})\otimes_{\bql}\bqlk$ est un facteur direct de $(R\Gm(X_{\g},\omega_{X_{\g}})\otimes_{\bql}\bqlk)\otimes_{\bql[X_*(T)]^{W}}\bql [X_*(T)]\cong (R\Gm(X_{\g},\omega_{X_{\g}})\otimes_{\bql}\bqlk)\otimes_{\bql}\bql[W]$.

Notons $B_{\g}=R\Gm(X_{\g},\omega_{X_{\g}})$, en tensorisant par $\otimes_{\bql[X_{*}(T)]}\bql$, le complexe $(R\Gm(X_{\g},\omega_{X_{\g}})\otimes_{\bql}\bqlk)\otimes_{\bql}\bql[W]$, il nous suffit donc de vérifier que:
\begin{equation}
(B_{\g}\otimes_{\bql}\bqlk)\otimes_{\bql[X_{*}(T)]^{W}}\bql\in\mathstrut^{p}\cD^{\geq-\nu_{w,r}}(\pt).
\label{gm-fib}
\end{equation}
 D'après \ref{act-sph}, la restriction de l'action de $X_*(T)$ à $\bql[X_*(T)]^{W}$ sur $B_{\g}$ se factorise par \eqref{sigma}, on obtient ainsi que:
\[(B_{\g}\otimes_{\bql}\bqlk)\otimes_{\bql[X_*(T)]^{W}}\bql\cong( B_{\g}\otimes_{\bql[\pi_{0}(P_{a})]}\bql[\pi_{0}(P_a)]\otimes_{\bql}\bqlk) \otimes_{\bql[X_*(T)]^{W}}\bql[\pi_{0}(P_a)]\otimes_{\bql[\pi_{0}(P_a)]}\bql.\]
On a alors le lemme clé suivant:
\bprop\label{p-calc}
Soit un caractère $\kappa:X_*(T)\ra\bZ/N\bZ$ pour $N\in\NN^*$. Considérons la DG-algèbre $\Kos(P_a)=\bigoplus\La^{i}(\bql[\pi_{0}(P_a)])[i]$ avec des différentielles nulles,  alors 
\[\bql[\pi_{0}(P_a)]\otimes_{\bql}\bqlk \otimes_{\bql[X_*(T)]^{W}}\bql[\pi_{0}(P_a)]\cong\bqlk\otimes_{\bql} \Kos(P_a),\]
si $\kappa_{\vert \bql[X_*(T)]^{W}}$ se factorise par $\pi_0(P_a)$ et 0 sinon.
\eprop
\bpf
Posons $A=\bql[X_*(T)]$, $A^{0}=\bql[X_*(T)]^{W}$ et $B=\bql[\pi_{0}(P_a)]$. On suppose que la restriction de $\kappa$ à $\bql[X_*(T)]^{W}$ ne se factorise pas par $\bql[\pi_{0}(P_a)]$.
Tout d'abord, en choisissant un relèvement $\sigma'_{\g}$ de $\sigma_{\g}$ tel que \eqref{w-lift}, il est équivalent  de dire que $\kappa$ ne se factorise par aucun de ces relèvements. Fixons un tel relèvement et notons $K=\Ker(\sigma'_{\g})\simeq\bZ^{l}$ pour un certain $l$ et on se sert de $\sigma'$ pour voir $B$ comme un $A$-module.
La suite exacte:
\[1\ra K\ra X_{*}(T)\stackrel{\sigma'_{\g}}{\rightarrow}\pi_{0}(P_a)\ra 1\]
fournit un isomorphisme de complexes:
\[\bqlk\otimes_{\bql[X_{*}(T)]}\bql[\pi_{0}(P_a)]\simeq\bqlk\otimes_{\bql[K]}\bql.\]
Par hypothèse $\kappa$ se restreint à $K$ en un caractère non-trivial. L'énoncé se déduit alors du lemme suivant:
\blem
Soit $\kappa:\bZ^l\ra\bZ/N\bZ$ un caractère non-trivial, alors le complexe des coinvariants $\bqlk\otimes_{\bql[\bZ^{l}]}\bql$ est nul.
\elem
\bpf
En effet, on voit $\bqlk$ comme le produit tensoriel extérieur de la représentation triviale sur $\bZ^{l-1}$ avec la représentation $\kappa$ de $\bZ$, de telle sorte que d'après Künneth, on a:
\[R\Gm(\bZ^{l},\bqlk)\cong R\Gm(\bZ^{l-1},\bql)\otimes_{\bql}R\Gm(\bZ,\bqlk).\]
Il suffit donc de voir que $R\Gm(\bZ,\bqlk)$ est nul.
On utilise alors le fait que le complexe
\[0\ra\bql[\bZ]\stackrel{t-1}{\rightarrow}\bql[\bZ]\ra 0\]
est une résolution du module trivial $\bql$, où $t$ est le générateur de $\bql[\bZ]$. Ainsi,  on en déduit que:
\[H_0(\bZ,\bqlk)=(\bqlk)_{\bZ},  H_{1}(\bZ,\bqlk)=(\bqlk)^{\bZ},\]
et les deux termes sont nuls comme $\kappa$ est non-trivial.
\epf
Si maintenant $\kappa$ se factorise par $\pi_0(P_a)$, alors comme $\bqlk\simeq\bqlk\otimes_{\bql[\pi_0(P_a)]}\bql[\pi_0(P_a)]$, il suffit de calculer $B\otimes_{A^0}B$.
On considère l'immersion fermée $\iota:\Omega_{a}=\Spec(B)\subset \check{T}/W=\Spec(A^0)$ qui est une immersion fermée entre schémas lisses, comme $G$ est simplement connexe, de codimension $c(\g)$ d'après \cite[Prop.3.9.2]{N}. En particulier, c'est une immersion régulière et le complexe de Koszul $\Kos(\iota)$ fournit une résolution libre de $B$ comme $A$-module d'après \cite[Tag. 062F]{Sta} de telle sorte que 
\[B\otimes^{L}_{A^{0}}B\cong\Kos(\iota)\otimes_{A^{0}}B\cong\bigoplus\limits_{i=0}^{c(\g)}\La^{i}B^{c(\g)}[i]\]
avec des différentielles nulles, ce qui conclut.
\epf
On peut maintenant terminer la preuve de \ref{coinv}:
\bpf

Il résulte de \ref{p-calc} que pour établir \eqref{gm-fib}, on peut supposer que $\kappa$ se factorise par $\pi_{0}(P_a)$.
D'après Kazhdan-Lusztig \cite[3.4]{N}, il existe alors un sous-groube libre de rang maximal $\La_{\g}\subset\pi_{0}(P_a)$ qui agit librement sur $X_{\g}$ de telle sorte que $Y_{\g}=X_{\g}/\La_{\g}$ est un espace algébrique propre et comme $\kappa$ est de torsion, on dispose d'un système local de rang $\cL_{\kappa}$ sur $Y_{\g}$.
On en déduit ainsi que:
\[(B_{\g}\otimes_{\bql}\bqlk)\otimes_{\bql[\La_{\g}]}\bql\cong R\Gm(Y_{\g},\cL_{\kappa})\]
et comme $\La_\g$ est de rang maximal, $(B_{\g}\otimes_{\bql}\bqlk)\otimes_{\bql[\pi_0(P_a)]}\bql$ est donc un facteur direct de $R\Gm(Y_{\g},\cL_{\kappa})$. En combinant avec \ref{p-calc}, on obtient donc que
\begin{equation}
(B_{\g}\otimes_{\bql}\bqlk)\otimes_{\bql[X_*(T)]^{W}}\bql\in\cD^{\geq-(2\delta(\g)+c(\g))}
\label{inegalité}
\end{equation}
et $\nu_{w,r}\geq 2\delta(\g)+c(\g)$ avec une inégalité stricte si $d(w,r)\geq 1$, ce qui conclut.
\epf
\brem\label{G-supp}
Il résulte de la dernière égalité que les strates qui peuvent donc intervenir comme supports sont les $\kC_{w,r}$ telles que $d(w,r)=0$ et il n'est pas clair que l'ensemble de ces strates forme un ouvert.
\erem
\subsection{Morphismes de spécialisation}
Le but de cette section est de construire des morphismes de spécialisation pour l'homologie des fibres de Springer affines. Cela permet en particulier de pouvoir relier l'homologie des fibres pour des éléments elliptiques et des éléments déployés.

Soit un $k$-schéma de type fini $X$ et $\eta, x\in X$ et $\bar{\eta}, \bar{x}$ des points géométriques au-dessus de $\eta$ et $x$, une spécialisation $\bar{\eta}\rightsquigarrow\bar{x}$ de codimension $r$, est un $X$-morphisme $\bar{\eta}\ra X_{(\bar{x})}$  dont l'adhérence de l'image de $\bar{y}$ est de dimension $r$.
Pour tout point géométrique $\bar{x}$ de $X$ d'image $x$, on note $d(\bar{x})=\dim\overline{\{x\}}$.
Pour tout faisceau $K\in\cD(X)$ et un point géométrique $\bar{x}$ de $X$, notons $R\Gm_{\bar{x}}K=\Gm(\bar{x},i_{\bar{x}}^{!}K)=\bD(K)_{\bar{x}}^{\vee}[2d(\bar{x})](d(\bar{x}))\in\Vect_{\bql}$, où ici $\Vect_{\bql}$ désigne la $\infty$-catégorie des complexes de $\bql$-espaces vectoriels et $\bD$ est le dual de Verdier sur $X$.

\begin{prop}\label{sp-morph}
Soit un $k$-schéma de type fini $X$  et une spécialisation de codimension $r$  $\bar{\eta}\rightsquigarrow\bar{x}$, soit $K\in\cD(X)$, alors on a un morphisme de spécialisation:
\begin{equation}
sp:R\Gm_{\bar{\eta}}K\ra R\Gm_{\bar{x}}K(r)[2r].
\label{sp-K}
\end{equation}
\end{prop}
\bpf
Les foncteurs étant continus, on se ramène au cas où $K\in\cD_c(X)$. Dans ce cas pour tout $L\in\cD_c(X)$, on a un morphisme de spécialisation au niveau des fibres:
\[L_{\bar{x}}\ra L_{\bar{\eta}}\]
et on prend $L=\bD(K)$ le dual de Verdier de $K$,  soit:
\[\bD(K)_{\bar{x}}\ra\bD(K)_{\bar{\eta}}\]
et on prend le dual de cette flèche pour obtenir:
\[R\Gm_{\bar{\eta}}K[-2d(\bar{\eta})](-d(\bar{\eta}))\cong(\bD(K)_{\bar{\eta}})^{\vee}\ra(i_{\bar{x}}^*\bD(K))^{\vee}.\]
d'où  comme $d(\bar{\eta})=r+d(\bar{x})$, une flèche:
\[sp:R\Gm_{\bar{\eta}}K\ra R\Gm_{\bar{x}}K(r)[2r].\]
\epf
Le morphisme \eqref{sp-K} est fonctoriel par rapport à tout morphisme étale $f:X'\ra X$ de telle sorte que pour tout groupe discret $\Gm$, en prenant $f:X\times\Gm\ra X$, on obtient que pour tout $K\in\cD^{\Gm}(X)$, le morphisme \eqref{sp-K} est $\Gm$-équivariant.

On obtient alors l'énoncé suivant, on le formule dans le cas de la droite affine, mais il vaut pour tout $k$-schéma de type fini.
\bprop\label{sp-sprin}
Considérons un morphisme $\phi:\ab^1\ra\kC$ alors pour tout morphisme de spécialisation  $\bar{\eta}\rightsquigarrow\bar{x}$ de $\ab^1$,
on a des flèches canoniques $\widetilde{W}$-équivariantes:
\begin{equation}
sp:R\Gm(X_{\g_{\eta}},\omega_{X_{\g_{\eta}}})\ra R\Gm(X_{\g_{s}},\omega_{X_{\g_{s}}})[2].
\label{sp-2}
\end{equation}
\eprop

\bpf
En effet, on forme le carré cartésien:
$$\xymatrix{X_{\la}\ar[d]_{f_{\la}}\ar[r]^{\ti{\phi}}&\wkC\ar[d]^{f}\\\ab^{1}\ar[r]^{\phi}&\kC}.$$
Ainsi, par changement de base propre \ref{Base}, on obtient :
\[\phi^{!}\cS=\ti{\phi}^{!}f_{!}\omega_{\wkC}\cong f_{\la,!}\ti{\phi}^!\omega_{\wkC}\cong f_{\la,!}\omega_{X_{\la}}\]
et on applique alors \ref{sp-morph} à $K=f_{\la,!}\omega_{X_{\la}}$ et à nouveau \ref{Base} pour obtenir une flèche  $\widetilde{W}$-équivariante:
\[sp:R\Gm(X_{\g_{\eta}},\omega_{X_{\g_{\eta}}})\ra R\Gm(X_{\g_{s}},\omega_{X_{\g_{s}}})[2].\]
\epf
\bexa
En utilisant l'action de $\widetilde{W}$, on peut obtenir des restrictions sur l'image de cette application.
Soit $G=PGL_2$, on considère une famille indexée par le paramètre $\la$:
\[a_{\la}:=X^2-(\la t^{2n}+t^{2n+1})\in\kc(k[\la][[t]]).\]
Si $\la\neq 0$, alors l'élément est déployé et $\pi_{0}(P_{a_{\la}})=\bZ$ et pour $\la=0$, on a $\pi_{0}(P_{a_{0}})=\{1\}$.
On a alors une flèche $\widetilde{W}$- équivariante :
\[R\Gm(X_{\g_{\eta}},\omega_{X_{\g_{\eta}}})\ra R\Gm(X_{\g_{s}},\omega_{X_{\g_{s}}})[2]\]
et en utilisant \ref{act-sph}, elle se factorise donc par les coinvariants:
\[R\Gm(X_{\g_{\eta}},\omega_{X_{\g_{\eta}}})_{\bZ}\ra R\Gm(X_{\g_{s}},\omega_{X_{\g_{s}}})[2].\]
\eexa

\brem
On s'attend à ce que le morphisme \eqref{sp-2} soit également $\cP$-équivariant, ce qui revient à considérer des faisceaux sur le champ quotient $[\ab^{1}/\cP]$. Malheureusement la construction de \ref{sp-morph} nécessite d'utiliser la dualité de Verdier, qui n'est pas disponible sur ce champ quotient.
\erem

\address{ 
  \bigskip
  \footnotesize

  (A.\ Bouthier) \textsc{IMJ-PRG, Sorbonne Université, 4 Place Jussieu, Paris, 75005 France}\par\nopagebreak
  \textit{E-mail address}: \texttt{alexis.bouthier@imj-prg.fr}}


\begin{thebibliography}{}
\bibitem{Alp}
J. Alper.
\newblock
Adequate moduli spaces and geometrically reductive group schemes.
\newblock
Algebraic Geometry 1 (4) pp. 489-531, (2014).


\bibitem{Bez}
R. Bezrukavnikov.
\newblock The dimension of the fixed points set on affine flag manifolds.
\newblock{\em Mathematical Research Letters} 3,  185-189, (1996).


\bibitem{BBD}
A. Beilinson, J. Bernstein, P. Deligne.
\newblock Faisceaux pervers.
\newblock Analyse et topologie sur les espaces singuliers, Astérisque vol. 100, SMF Paris, 1982.

\bibitem{BeKV}
R. Bezrukavnikov, D. Kazhdan, Y. Varshavsky.
\newblock A categorical approach to the stable
center conjecture.
\newblock{\em Astérisque} 369, pp. 27-97, 2015 SMF, Paris.

\bibitem{Bt}
A. Bouthier.
\newblock Support singulier et homologie des fibres de Springer affines, https://arxiv.org/abs/arXiv:2202.12017.

\bibitem{BC}
A. Bouthier, K. \v{C}esnavi\v{c}ius.
\newblock Torsors on loop groups and the Hitchin fibration.
\newblock{\em Annales scientifiques de l'ENS} 55, no. 3, pp. 791-864,  2022.


\bibitem{BKV}
A. Bouthier, D. Kazhdan, Y. Varshavsky.
\newblock
Perverse sheaves on infinite-dimensional stacks, and affine Springer Theory.
\newblock{\em Advances in Maths}, vol. 408, part. A, pp., Oct. 2022.

\bibitem{CG}
N. Chriss and V. Ginzburg.
\newblock
 Representation theory and complex geometry, Birkhäuser, Boston, 1997.

\bibitem{SGA4.5}
P. Deligne.
\newblock Cohomologie étale SGA 4 1/2.
\newblock Lect. Notes in Math. 569, Springer-Verlag, pp. 233-251, 1977.

\bibitem{DG}
M. Demazure, P. Gabriel.
\newblock Groupes algébriques. Tome I: Géométrie algébrique, généralités, groupes
commutatifs.
\newblock Avec un appendice Corps de classes local par Michiel Hazewinkel. Masson Editeur, Paris, North-Holland Publishing Co., Amsterdam (1970).

\bibitem{Dr}
V. Drinfeld.
\newblock
Infinite-dimensional vector bundles in algebraic geometry: an introduction.
\newblock The Unity of Mathematics, pp. 263-304, Birkhäuser, Boston 2006.


\bibitem{Gai}
D. Gaitsgory.
\newblock Space of rational maps.
\newblock{\em Inventiones Math.}, vol. 191, pp. 91-196, 2013.

\bibitem{GR}
D. Gaitsgory, N. Rozenblyum.
\newblock
DG-indschemes.
\newblock Contemporary Mathematics 610, pp. 139-251, 2014.

\bibitem{GKM}
M. Goresky, R. Kottwitz, R. McPherson.
\newblock
Codimension of root valuation strata.
\newblock Pure and Applied Mathematics Quarterly 5 , 1253-1310, 2009.

\bibitem{GKM2}
M. Goresky, R. Kottwitz, R. McPherson.
\newblock Purity of equivalued affine Springer fibers.
\newblock{\em Representation Theory}, vol. 10, pp. 130-146, (2006).

\bibitem{HR}
T. Haines, T. Richarz.
\newblock
The Test Function Conjecture for Local Models of Weil-restricted groups.
\newblock{\em Compositio Mathematica} 156, 1348-1404, 2020.


\bibitem{EGAIV}
A. Grothendieck, J. Dieudonné.
\newblock Éléments de géométrie algébrique. IV: Étude locale des schémas et des morphismes de schémas
Quatrième partie.
\newblock{\em Publ. Math. IHES}, vol. 32, (1967).

\bibitem{SGA1}
A. Grothendieck.
\newblock Séminaire de Géométrie Algébrique du Bois Marie: Revêtements étales et groupe
fondamental.
\newblock Lecture notes in mathematics 224, Springer-Verlag (1971).

\bibitem{SGA3}
P. Gille and P. Polo.
\newblock
 Schémas en groupes (SGA 3). Tome I. Propriétés générales des
schémas en groupes.
\newblock
 Documents Mathématiques (Paris), 7, SMF, Paris, 2011 . Séminaire de Géométrie Algébrique du Bois Marie
1962-64 dirigé par M. Demazure and A. Grothendieck.

\bibitem{KL}
D. Kazhdan, G. Lusztig.
\newblock Fixed point varieties on affine flag manifolds.
\newblock{\em Israël J. Math. 62}, no. 2, pp. 129-168, (1988).

\bibitem{KW}
R. Kiehl, R. Weissauer.
\newblock Weil conjectures, Perverse sheaves and $\ell$-adic Fourier transform.
\newblock Ergebnisse der Mathematik und ihrer Grenzgebiete, vol.42, Springer, Berlin.

\bibitem{Sta}
A. J. de Jong et al.
\newblock The Stacks Project. Available at http://stacks.math.columbia.edu.

\bibitem{LN}
G. Laumon, B.C. Ngô.
\newblock
Le lemme fondamental pour les groupes unitaires.
\newblock{\em Ann. of Math.}, vol. 168, pp. 477-573, 2008.

\bibitem{LZ1}
Y.Liu, W. Zheng.
\newblock Enhanced six operations and base change theorem for sheaves on Artin
stacks. preprint, arXiv:1211.5948.


\bibitem{LZ2}
Y.Liu, W. Zheng.
\newblock
Enhanced adic formalism and perverse t-structures for higher Artin stacks.
preprint, arXiv:1404.1128.


\bibitem{Lu1}
J. Lurie.
\newblock
Higher Topos theory.
\newblock Annals of Mathematics Studies, vol. 170, Princeton University Press, Princeton, NJ, 2009.


\bibitem{Lu2}
J. Lurie.
\newblock
Higher algebra. 
\newblock disponible sur https://www.math.ias.edu/~lurie/papers/HA.pdf.


\bibitem{Lu}
G. Lusztig.
\newblock Affine Weyl groups and conjugacy classes in Weyl groups.
\newblock{\em Transform. Groups}, no. 1-2, pp. 83-97, 1996.

\bibitem{Mi}
S. A. Mitchell.
\newblock
Quillen's theorem on buildings and the loops on a symmetric space.
\newblock{\em Enseign. Math. }(2) 34 (1-2), pp. 123-166, (1988).

\bibitem{N}
B.C. Ngô.
\newblock Le lemme fondamental pour les algèbres de Lie.
\newblock{\em Publ. Math. Inst. Hautes Études Sci. 111}, pp. 1-169, 2010.

\bibitem{N2}
B.C. Ngô.
\newblock Endoscopie et Fibration de Hitchin.
\newblock{\em Invent. Math.}, pp. 399-453, vol. 164, 2006.


\bibitem{Roz}
N. Rozenblyum.
\newblock
Filtered colimits of $\infty$-categories,
\newblock disponible sur http://www.math.harvard.edu/~gaitsgde/GL/colimits.pdf.


\bibitem{Tho}
R. Thomason.
\newblock Equivariant resolution, Linearization and Hilbert's Fourteenth Problem over arbitrary base schemes.
\newblock{\em Advances in Maths}, vol. 65, 16-34 (1987).


\bibitem{Yun}
Z. Yun.
\newblock Global Springer theory.
\newblock Adv. Math. 228  no. 1, 266-328, (2011).


\bibitem{YunII}
Z. Yun.
\newblock The spherical part of the local and global Springer actions.
\newblock Math. Ann. 359 (2014), no. 3-4, pp. 557-594.

\end{thebibliography}
\end{document}